\documentclass[11pt]{article}

\newcommand{\documentdate}{21 XI 2021}
\usepackage{latexsym,varioref,amsmath,mathtools}

\newcommand{\numsection}[1]{\section{#1}\setcounter{equation}{0}}

\newcommand*\varybar[1]{%
   \hbox{%
     \vbox{%
       \hrule height 0.5pt 
       \kern0.5ex
       \hbox{%
         \kern-.2em
         \ensuremath{#1}%
         \kern-0.0em
       }%
     }%
   }%
} 
\newcommand*\barDT{%
   \kern0.2em
   \hbox{%
     \vbox{%
       \hrule height 0.5pt 
       \kern0.5ex
       \hbox{%
         \kern-0.2em
         \ensuremath{\Delta T}%
         \kern-0.0em
       }%
     }%
   \kern0.0em
   }%
} 
\newcommand*\barphi{%
   \kern0.1em
   \hbox{%
     \vbox{%
       \hrule height 0.5pt 
       \kern0.5ex
       \hbox{%
         \kern-0.1em
         \ensuremath{\phi}%
         \kern-0.0em
       }%
     }%
   }%
   \kern0.0em
} 
\newcommand*\barf{%
   \kern0.2em
   \hbox{%
     \vbox{%
       \hrule height 0.5pt 
       \kern0.5ex
       \hbox{%
         \kern-0.2em
         \ensuremath{f}%
         \kern-0.0em
       }%
     }%
   }%
   \kern0.0em
} 

\newcommand*\barnablaxif{%
   \kern0.2em
   \hbox{%
     \vbox{%
       \hrule height 0.3pt 
       \kern0.3ex
       \hbox{%
         \kern-0.0em
         \ensuremath{\nabla_x^i f}%
         \kern-0.0em
       }%
     }%
   }%
   \kern0.0em
} 
\newcommand*\barnablaxjf{%
   \kern0.2em
   \hbox{%
     \vbox{%
       \hrule height 0.3pt 
       \kern0.3ex
       \hbox{%
         \kern-0.0em
         \ensuremath{\nabla_x^j f}%
         \kern-0.0em
       }%
     }%
   }%
   \kern0.0em
} 
\newcommand{\al}[1]{{\footnotesize{\sf #1}}}
\newcommand{\tal}[1]{{\normalsize {\sf #1}}}
\newcommand{\DT}{\Delta T}
\newcommand{\accuracy}{{\tt accuracy}}

\newcommand{\sufficient}{{\tt sufficient}}
\newcommand{\insufficient}{{\tt insufficient}}
\newcommand{\absolute}{{\tt absolute}}
\newcommand{\relative}{{\tt relative}}

\newcommand{\dsk}{\delta_{s_k}}
\newcommand{\dskj}{\delta_{s_k,j}}

\newcommand{\dskl}{\delta_{s_k,\ell}}

\newcommand{\barm}{\overline{m}  }  
\newcommand{\aacc}{\varphi}         
\newcommand{\countacc}{k_\aacc}      
\newcommand{\kapsteptwo}{\kappa_{\mbox{\scriptsize step2}}} 
\newcommand{\epsmin}{\epsilon_{\min}}
\newcommand{\s}[1]{^{\mbox{\protect\tiny #1}}}
\newcommand{\sub}[1]{_{\mbox{\protect\tiny #1}}}
\newcommand{\kat}[1]{\kap{\sf #1}}
\newcommand{\kats}[1]{\kat{#1}\s{S}}
\newcommand{\kata}[1]{\kat{#1}\s{A}}

\newcommand{\katc}[1]{\kat{#1}\s{C}}

\newcommand{\kate}[1]{\kat{#1}\s{E}}
\newcommand{\katf}[1]{\kat{#1}\s{F}}
\newcommand{\evbndf}[1]{\mbox{\textrm{\textup{\small\texttt{\#}}}}_{\mbox{\tiny {\sf #1}}}\s{F}}
\newcommand{\evbndd}[1]{\mbox{\textrm{\textup{\small\texttt{\#}}}}_{\mbox{\tiny {\sf #1}}}\s{D}}
\newcommand{\flow}{f\sub{low}}
\newcommand{\beqn}[1]{\begin{equation}\label{#1}}
\newcommand{\eeqn}{\end{equation}}
\renewcommand{\Re}{\hbox{I\hskip -2pt R}}
\newcommand{\smallRe}{\hbox{\footnotesize I\hskip -2pt R}}
\newcommand{\eqdef}{\stackrel{\rm def}{=}}
\newcommand{\bigfrac}[2]{\frac{\displaystyle #1}{\displaystyle #2}}
\newcommand{\bigsum}{\displaystyle \sum}
\newcommand{\bigmax}{\displaystyle \max}
\newcommand{\bigmin}{\displaystyle \min}
\newcommand{\req}[1]{(\ref{#1})}
\newcommand{\tim}[1]{\;\; \mbox{#1} \;\;}
\newcommand{\ii}[1]{\{1, \ldots, #1 \}}
\newcommand{\iiz}[1]{\{0, \ldots, #1 \}}
\newcommand{\iibe}[2]{\{ #1, \ldots, #2 \}}
\newcommand{\barL}{\overline{L}}
\newcommand{\barT}{\overline{T}}
\newcommand{\calA}{{\cal A}}
\newcommand{\calB}{{\cal B}}
\newcommand{\calE}{{\cal E}}
\newcommand{\calO}{{\cal O}}
\newcommand{\calS}{{\cal S}}
\newcommand{\calT}{{\cal T}}
\newcommand{\calU}{{\cal U}}
\newcommand{\eo}[1]{\calO\left(#1\right)}
\newcommand{\kap}[1]{\kappa_{\mbox{\rm \tiny #1}}}

\newcounter{algo}[section]
\renewcommand{\thealgo}{\thesection.\arabic{algo}}

\newcommand{\algo}[3]{\refstepcounter{algo}
\begin{center}\begin{figure}[htbp]
\framebox[\textwidth]{
\parbox{0.95\textwidth} {\vspace{\topsep}
{\bf Algorithm \thealgo : #2}\label{#1}\\
\vspace*{-\topsep} \mbox{ }\\
{#3} \vspace{\topsep} }}
\end{figure}\end{center}}

\newcommand{\balgo}[4]{\refstepcounter{algo}
\begin{center}\begin{figure}[htbp]
\vspace{\baselineskip}\noindent\hbox{%
  \lower\fboxrule\hbox{\vbox{\hrule\hbox{\vrule \kern-\fboxrule \vbox{%
  \vspace{\topsep} \noindent\hspace{2\fboxsep}\parbox{\thmw}{\vspace{0.5\topsep}
  {\bf Algorithm \thealgo : #2}\label{#1}\\
  \vspace*{-\topsep} \mbox{ }\\
  {\rm #3}\vspace{-\lastskip}}
  \hspace{\fboxsep}}\kern-\fboxrule \vrule }}}}
\end{figure}\end{center}
\newpage \hbox{%
  \lower\fboxrule\hbox{\vbox{\hbox{\vrule \kern-\fboxrule \vbox{%
  \noindent\hspace{2\fboxsep}\parbox{\thmw}{\rm #4}\hspace{\fboxsep}
  \vspace{\fboxsep}}\kern-\fboxrule \vrule }\hrule }}}\vspace{\baselineskip}
}

\newcommand{\balgob}[5]{\refstepcounter{algo}
\begin{center}\begin{figure}[htbp]
\vspace{\baselineskip}\noindent\hbox{%
  \lower\fboxrule\hbox{\vbox{\hrule\hbox{\vrule \kern-\fboxrule \vbox{%
  \vspace{\topsep} \noindent\hspace{2\fboxsep}\parbox{\thmw}{\vspace{0.5\topsep}
  {\bf Algorithm \thealgo : #2}\label{#1}\\
  \vspace*{-\topsep} \mbox{ }\\
  {\rm #3}\vspace{-\lastskip}}
  \hspace{\fboxsep}}\kern-\fboxrule \vrule }}}}
\end{figure}\end{center}
\newpage \hbox{%
  \lower\fboxrule\hbox{\vbox{\hbox{\vrule \kern-\fboxrule \vbox{%
  \noindent\hspace{2\fboxsep}\parbox{\thmw}{\rm #4}\hspace{\fboxsep}
  \vspace{\fboxsep}}\kern-\fboxrule \vrule } }}}\vspace{\baselineskip}
\newpage \hbox{%
  \lower\fboxrule\hbox{\vbox{\hbox{\vrule \kern-\fboxrule \vbox{%
  \noindent\hspace{2\fboxsep}\parbox{\thmw}{\rm #5}\hspace{\fboxsep}
  \vspace{\fboxsep}}\kern-\fboxrule \vrule }\hrule }}}\vspace{\baselineskip}
}

\newtheorem{theorem}{Theorem}[section]
\newtheorem{lemma}[theorem]{Lemma}
\newcommand{\llem}[2]{\vspace{\baselineskip} 
\noindent\framebox[\textwidth]{\parbox{0.95\textwidth}{
    \begin{lemma} \label{#1} \rm #2 \end{lemma} } } \vspace{\baselineskip} }

\newcommand{\lbthm}[3]{\vspace{\baselineskip}\noindent\hbox{%
  \lower\fboxrule\hbox{\vbox{\hrule\hbox{\vrule \kern-\fboxrule \vbox{%
  \vspace{\fboxsep} \noindent\hspace{2\fboxsep}\parbox{\thmw}{
  \begin{theorem}\label{#1}{\rm #2}\end{theorem}\vspace{-\lastskip}}
  \hspace{\fboxsep}}\kern-\fboxrule \vrule }}}}\newpage \hbox{%
  \lower\fboxrule\hbox{\vbox{\hbox{\vrule \kern-\fboxrule \vbox{%
  \noindent\hspace{2\fboxsep}\parbox{\thmw}{\rm #3}\hspace{\fboxsep}
  \vspace{4\fboxsep}}\kern-\fboxrule \vrule }\hrule }}}\vspace{\baselineskip}
}

\newcommand{\bpr}{{\bf Proof.} \hspace{1.5mm}}
\newcommand{\epr}{\hfill $\Box$ \vspace*{1em}}
\newcommand{\proof}[1]{
\begin{list}{}{
\setlength{\topsep}{0.0pt}
\setlength{\partopsep}{0.0pt}
\setlength{\leftmargin}{0.025\textwidth}
\setlength{\rightmargin}{0.5\leftmargin}
\setlength{\labelwidth}{0.5\leftmargin}
\setlength{\labelsep}{0.25\leftmargin}}
\item \bpr #1 \epr \noindent
\end{list}}
\newlength{\thmw}
\newcommand{\arr}[2]{\begin{array}{#1}#2\end{array}}
\newcommand{\sfrac}[2]{{\scriptstyle \frac{#1}{#2}}}
\newcommand{\half}{\sfrac{1}{2}}
\newcommand{\quarter}{\sfrac{1}{4}}
\newcommand{\ms}{\;\;\;\;}

\title{An adaptive regularization algorithm for unconstrained optimization
  with inexact function and derivatives values}

\author{
N. I. M. Gould\thanks{Computational Mathematics Group,
   STFC-Rutherford Appleton Laboratory,
   Chilton OX11 0QX, England. Email:  nick.gould@stfc.ac.uk .
   The work of this author was supported by EPSRC grant EP/M025179/1}
~and~Ph. L. Toint\thanks{Namur Center for Complex Systems (naXys),
   University of Namur, 61, rue de Bruxelles, B-5000 Namur, Belgium.
   Email: philippe.toint@unamur.be}
}

\date{\documentdate}

\begin{document}


\maketitle

\begin{abstract}
  An adaptive regularization algorithm for unconstrained nonconvex
  optimization is proposed that is capable of handling inexact
  objective-function and derivative values, and also of providing approximate
  minimizer of arbitrary order. In comparison with a similar algorithm 
  proposed in Cartis, Gould, Toint (2021), its distinguishing feature is that 
  it is based on controlling the relative error between the model and 
  objective values. A sharp
  evaluation complexity complexity bound is derived for the new algorithm.
\end{abstract}

\vspace*{0.5cm}
\textbf{Keywords: } nonconvex optimization, inexact functions and
derivative values, evaluation complexity, adaptive regularization.

  \numsection{Introduction: motivation, context, definitions}

We consider the unconstrained minimization problem
\beqn{problem}
\min_{x\in\smallRe^n} f(x),
\eeqn
where $f$ is a $p$-times continuously differentiable function from $\Re^n$ to 
$\Re$. In practice it sometimes happens that the value of the objective
function $f$ and/or those of its derivatives may only be computed inexactly,
for a variety of reasons. It might be that the
derivatives are not available, and are estimated using, for example, 
finite differences. Or perhaps because the evaluation is subject to 
intentional noise.
For example, the values in question might be computed by some kind of 
experimental process whose accuracy can be adjusted, with the understanding 
that more accurate values may be, sometimes substantially, more expensive 
in terms of computational effort.  A related case is when objective-function or
derivative values result from some (one hopes, convergent) 
iteration---obtaining more accuracy is possible by letting the iteration 
converge further, but again at the price of possibly significant additional 
computing.  A third possibility, currently much in vogue in the context of 
machine learning, is when the values of the objective function and/or 
its derivatives are obtained by sampling---say from among the terms 
of a sum of functions involving an enormous number of
them. Again, using a larger sample size results in probabilistically better
accuracy, but at a cost. In this report, we are particularly interested in a
fourth case of growing relevance in high-performance computing, which is the
emerging field of ``variable accuracy'' or ``multi-precision'' optimization 
\cite{Baboetal09,Leyfetal16,BellGratRicc18,GratToin20},
where the computation of the objective function and its derivative(s) are
intentionally truncated to significantly fewer digits than would be required
for an accurate calculation. Doing so allows the use of specialised processing
units whose chip surface, and hence power consumption, is much less than is
needed for standard double precision arithmetic \cite{GalaHoro11}.

It is therefore important to design and analyze algorithms which are tolerant
to inexactness or noise.  This subject is not new.  The main ideas developed
in this report have their origin in the paper \cite{BellGuriMoriToin19}, that
considers the complexity of finding weak approximate minimizers using an
algorithm similar, at least in spirit, to the algorithm discussed here. Our new
presentation and its associated analysis
merge the elaborate approximation techniques described in that paper with
techniques of \cite[Chapter~12]{CartGoulToin21} for computing (strong)
approximate minimizers and has been explored at length in
\cite[Chapter~13]{CartGoulToin21}. It also avoids using a dynamic relative
accuracy threshold present in \cite{BellGuriMoriToin19}. First-order
trust-region methods with inexact evaluations and explicit dynamic accuracy
have been described in \cite{Cart93} and in Section~10.6 of
\cite{ConnGoulToin00}. The complexity analysis for the first-order
case was also discussed in \cite{GratToin20}. 
 
The purpose of the report in hand is to present an alternative to one of the
algorithms described in \cite[Chapter~13]{CartGoulToin21}, namely
\al{AR$qp$EDA}, the adaptive regularisation
algorithm with explicit dynamic accuracy, and to provide a sharp
upper bound on its evaluation complexity. The variant we will consider here,
which we call \al{AR$qp$EDA2}, uses the same Explicit Dynamic Accuracy (EDA)
framework as \al{AR$qp$EDA}, but enforces its associated controls in a
different way.  Obviously, these sentences require some clarification as to
what these terms convey.

By ``explicit dynamic accuracy'', we mean that, during the optimization process,
\emph{the required values} (objective-function or derivatives) \emph{can
always be computed with an accuracy that is explicitly\footnote{See
\cite[Chap.~13]{CartGoulToin21} for the description and analysis of an
\emph{implicit} dynamic accuracy framework.} specified, before the
calculation, by the algorithm itself.}  It is also 
understood in what follows that \emph{the algorithm should require high
accuracy only if necessary}, but nonetheless \emph{guarantee final 
results to full requested accuracy}. In this situation, it is hoped that many 
function or derivative evaluations can be carried out with a fairly loose
accuracy (we will refer to these as ``inexact values''), thereby resulting in
a significantly cheaper optimization process.

Different kinds of optimization method can be designed to work within this
framework. We focus here an ``adaptive regularization'' algorithms(see
\cite{Grie81,Nest06,CartGoulToin11} among many other contributions), where, at
a given iteration, a regularized local Taylor model of the objective function
is minimized to define the next trial iterate.  The strength of the
regularization is then adaptively updated to guarantee convergence to
approximate minimizers. In addition, we will consider a variant of the
adaptive regularization technique whose purpose is to find such \emph{approximate
minimizers of arbitrary}, but given, \emph{order}. If an approximate
minimizer of order $q \ge 1$ is sought, this requires the Taylor model of the
objective function to have degree $p\geq q$. Such a model is of the form
\beqn{taylex}
T_{f,j}(x,s) \eqdef f(x) + \bigsum_{i=1}^j \frac{1}{i!}\nabla_x^i f(x)[s]^i
\equiv T_{f,j}(x,0) +  \bigsum_{i=1}^j \frac{1}{i!} [\nabla_v ^i T_{f,j}(x,v)]_{v=0}[s]^i
\eeqn
for perturbations $s$ around $x$, where $\nabla^j_xf(x)$ is the $j$-th order
tensor giving the the $j$-th derivative of $f$ at $x$, and where the notation
$\nabla^j_xf(x)[s]^j$ means that it is applied on $j$ copies of $s$

Given a Taylor model \req{taylex}, we also have to define what we
mean by an $\epsilon$-approximate minimizer of order $q$.  Following
\cite{CartGoulToin20a} (and \cite{CartGoulToin21}), we say that $x$ is 
an $\epsilon$-approximate minimizer of order $q$ whenever, for some
$\delta\in (0,1]^q$ and $\epsilon \in (0,1]^q$,
\beqn{optex}
\phi_{f,j}^{\delta_j}(x) \leq \epsilon_j \frac{\delta^j}{j!}
\tim{ for } j\in\ii{q},
\eeqn
where
\beqn{phiex}
\phi_{f,j}^{\delta_j}(x) = f(x) - \min_{d\in\smallRe^n,\,\|d\|\leq \delta_j}T_{f,j}(x,d).
\eeqn
In other words, this is the case when the none of the Taylor approximations
$T_{f,j}$ ($j\leq q$) at $x$ can be decreased by more than a scaled multiple
of $\epsilon_j$ (the right-hand side of \req{optex}) in a neighbourhood of $x$
of radius $\delta_j$. (The $\delta_j$ are called the \emph{optimality radii}.)
We again refer the reader to \cite{CartGoulToin20a} or \cite{CartGoulToin21}
for discussion of this optimality concept and of why it is a suitable
genralization of the standard low order optimality conditions to arbitrary
order.

In our new EDA framework, we have to be content with an inexact equivalent of
\req{taylex} given by
\beqn{ipa-barT-def}
\barT_{f,j}(x,s)
= \barf(x) + \bigsum_{i=1}^j \frac{1}{i!}\overline{\nabla_x^i f}(x)[s]^i
\equiv \barT_{f,j}(x,0)+\bigsum_{i=1}^j\frac{1}{i!}[\nabla_v ^i\barT_{f,j}(x,v)]_{v=0}[s]^i.
\eeqn
where, here and hereafter, we denote inexact quantities and approximations with an overbar.
It is therefore pertinent to investigate the effect of inexact derivatives on
\req{ipa-barT-def} and its uses. More specifically, we will be concerned with the  
Taylor {\em decrement} $\DT_{f,j}(x,s)$ at $x$ and for a step $s$, defined as
\beqn{ipa-taylor-dec}
\DT_{f,j}(x,s)
\eqdef  T_{f,j}(x,0) - T_{f,j}(x,s) 
= - \bigsum_{i=1}^j \frac{1}{i!} [\nabla_v ^i T_{f,j}(x,v)]_{v=0}[s]^i 
\equiv - \bigsum_{i=1}^j \frac{1}{i!}\nabla_x^i f(x)[s]^i.
\eeqn
While our traditional algorithms depend on this quantity, it is of course 
out of the question to use it in the present context, as we only have 
approximate values. But an obvious alternative is instead to consider 
the {\em inexact} Taylor decrement
\beqn{ipa-barDT}
\barDT_{f,j}(x,s)
\eqdef  \barT_{f,j}(x,0)-\barT_{f,j}(x,s)
= - \bigsum_{i=1}^j \frac{1}{i!} [\nabla_v^i \barT_{f,j}(x,v)]_{v=0}[s]^i,
\eeqn
itself resulting in an inexact version of \req{phiex} given by
\beqn{barphi-def}
\barphi_{f,j}^{\delta_{k,j}}(x) = - \max_{d\in\smallRe^n,\,\|d\|\leq \delta_j}\barDT_{f,j}(x,d).
\eeqn
In the ``Explicit Dynamic Accuracy'' (EDA) framework, we assume that the
conditions
\beqn{eda-vareps-j}
\|\overline{\nabla_x^if}(x)-\nabla_x^if(x)\| \leq \aacc_i
\eeqn
for degrees $i \geq 0$ of interest \emph{are enforced} when demanded by the 
algorithm, where $\aacc_i$ is the required absolute error bound on the
$i$-th derivative. It then follows that the error between the exact and
inexact Taylor expansions satisfies the bound 
\beqn{T-err0}
|\barDT_{f,j}(x,s)-\DT_{f,j}(x,s)|
\leq \bigsum_{i=1}^j \|\overline{\nabla_x^i f}(x)- \nabla_x^i f(x)\|\,
 \frac{\|s\|^{i}}{i!}\\*[2ex]
\leq \bigsum_{i=1}^j \aacc_{k,i} \frac{\|s\|^i}{i!},
\eeqn
using the triangle inequality and the requirement \req{eda-vareps-j}.

The distingushing feature of the \al{AR$qp$EDA2} algorithm---compared with
\al{AR$qp$EDA} of \cite[Chapter~13]{CartGoulToin21}---is that it enforces
convergence by controlling a (scaling independent) relative error bound  
\beqn{eda-barDT-acc}
|\barDT_{f,p}(x,s)-\DT_{f,p}(x,s)| \leq \omega \barDT_{f,p}(x,s),
\eeqn
for some fixed relative accuracy parameter $\omega \in (0,1)$. Most of the
difficulty in the forthcoming analysis results from the need to impose this
relative error bound and it is not obvious at this point how it can be
enforced using the absolute bounds \req{eda-vareps-j}. We now consider how
this can be achieved.

\numsection{Enforcing accuracy of the Taylor decrements}

For clarity, we temporarily neglect the iteration index $k$.
While there may be circumstances in which \req{eda-barDT-acc} can be enforced
directly, we consider here that the only control the user has on the accuracy
of $\barDT_{f,j}(x,s)$ is by imposing the bounds \req{eda-vareps-j}
on the \emph{absolute} errors of the derivative 
tensors $\{\nabla_x^if(x)\}_{i=1}^j$. As one may anticipate by examining
\req{eda-barDT-acc}, a suitable relative accuracy requirement can be
achieved so long as $\barDT_{f,j}(x,s)$ remains safely away from zero.
However, if exact computations are to be avoided, we may have to 
accept a simpler absolute accuracy guarantee when $\barDT_{f,j}(x,s)$ is 
small, but one that still guarantees our final optimality conditions. 

Of course, not all derivatives need to be inexact.  If 
derivatives of order $i \in \calE\subseteq \ii{q}$ are exact, then the
left-hand side of \req{eda-vareps-j} vanishes for $i\in \calE$ and
the choice $\aacc_i=0$ for $i\in \calE$ is perfectly adequate. 
However, we avoid the notational complication that making this
distinction would entail.

We start by describing a crucial tool that we use to achieve
\req{eda-barDT-acc}. This is the \al{CHECK} algorithm, stated as
Algorithm~\ref{eda-verify} \vpageref[below]{eda-verify}.
We use this to assess the relative model accuracy whenever needed in the
algorithms we describe later in this section.

To put our exposition in a general context, we suppose that we 
have an $r$-th degree Taylor series $T_r(x,v)$ of a given
function about $x$ in the direction $v$, along with an 
approximation $\barT_r(x,v)$ and its decrement $\barDT_r(x,v)$.  
Additionally, we suppose that a bound $\delta \geq \|v\|$ is given, and that
\emph{required} relative and absolute accuracies $\omega$ and $\xi>0$ 
are on hand; the relative accuracy constant $\omega \in (0,1)$ will 
fixed throughout the forthcoming algorithms, and we assume that it is 
given when needed in \al{CHECK}.
Finally, we assume that the \emph{current}
upper bounds $\{\aacc_j\}_{j=1}^r$ on absolute accuracies of the derivatives of
$\overline{T}_r(x,v)$ with respect to $v$ at $v=0$ are provided.
Because it will always be the case when we need it, 
for simplicity we will assume that $\barDT_r(x,v) \geq 0$.  

\algo{eda-verify}{Verify the accuracy of $\boldmath{\barDT_r(x,v)}$ 
(\al{CHECK})}
{
\[
\accuracy
= \mbox{\al{CHECK}}\Big(\delta,\barDT_r(x,v),\{\aacc_i\}_{i=1}^{r},\xi\Big).
\]
\begin{description}
\item[\hspace*{6.5mm}] 
  If  \vspace*{-3mm}
  \beqn{extDA-verif-term-2}
  \barDT_r(x,v) > 0
  \tim{ and }
  \sum_{i=1}^r \aacc_i \frac{\delta^i}{i!} \leq \omega \barDT_r(x,v),
  \vspace*{-3mm}
  \eeqn
  set \accuracy\ to \relative.
\item[\hspace*{6.5mm}] 
  Otherwise, if
  \vspace*{-3mm}
  \beqn{extDA-verif-term-3}
  \sum_{i=1}^r \aacc_i\frac{\delta^i}{i!} \leq \omega \xi \frac{\delta^r}{r!},
  \vspace*{-3mm}
  \eeqn
  set \accuracy\ to \absolute.
\item[\hspace*{6.5mm}] 
 Otherwise  set \accuracy\ to \insufficient.
  \end{description}
}

\newpage
\noindent
It will be convenient to say informally that 
\accuracy\ is \sufficient, if it is either \absolute\ or \relative.
We may formalise the accuracy guarantees that result from applying the
\al{CHECK} algorithm as follows.

\llem{eda-verify-l}{
Let $\omega \in (0,1]$ and $\delta, \xi$ and $\{\aacc_i\}_{i=1}^{r}>0$. 
Suppose that $\barDT_r(x,v) \geq 0$, that
\[\accuracy
= \mbox{\al{CHECK}}\Big(\delta,\barDT_r(x,v),\{\aacc_i\}_{i=1}^{r},\xi\Big),
\]
and that
\beqn{extDA-rvareps-j}
\Big\|\Big[\nabla_v^i \barT_r(x,v)\Big]_{v=0}
-\Big [\nabla_v ^i T_r(x,v)\Big]_{v=0}\Big\|
\leq \aacc_i \tim{for} i \in \ii{r}.
\eeqn
Then

\noindent
{\em (i)} \accuracy\ is \sufficient\ whenever
  \vspace*{-1mm}
  \beqn{extDA-max-abs-acc}
    \sum_{i=1}^r\aacc_i\frac{\delta^i}{i!} \leq \omega \xi \frac{\delta^r}{r!}.
  \vspace*{-2mm}
  \eeqn
 
\noindent
{\em (ii)} if \accuracy\ is \absolute, 
  \beqn{extDA-verif-prop-1}
  \max\Big[\barDT_r(x,w),
    \left|\barDT_r(x,w)- \DT_r(x,w)\right|\Big]
  \leq \xi \frac{\delta^r}{r!}
  \eeqn
  for all $w$ with $\|w\|\leq \delta$.

\noindent
{\em (iii)} if \accuracy\ is \relative, 
  $\barDT_r(x,v)>0$ and
  \beqn{extDA-verif-prop-2}
  \left|\barDT_r(x,w)- \DT_r(x,w)\right|
  \leq \omega \barDT_r(x,w),
  \tim{for all $w$ with $\|w\|\leq \delta$.} 
  \eeqn

}

\proof{
We first prove proposition \textup{(i)}, and assume that 
\req{extDA-max-abs-acc} holds, which clearly ensures that \req{extDA-verif-term-3} is satisfied. 
Thus either \req{extDA-verif-term-2} or \req{extDA-verif-term-3} must hold and
termination occurs, proving the first proposition.

It follows by definition of the decrements \req{ipa-taylor-dec} 
and \req{ipa-barDT}, the triangle inequality and \req{extDA-rvareps-j} that
\beqn{extDA-difft}
\arr{rl}{
\left|\barDT_r(x,w)- \DT_r(x,w)\right| 
& = \left| \bigsum_{i=1}^r \bigfrac{\left(
\left[\nabla_v^i \barT_r(x,v)\right]_{v=0}
-\left[\nabla_v ^i T_r(x,v)\right]_{v=0}
\right)[w]^i}{i!} \right| \\
& \leq \bigsum_{i=1}^r \bigfrac{\left\|
\left[\nabla_v^i \barT_r(x,v)\right]_{v=0}
-\left[\nabla_v ^i T_r(x,v)\right]_{v=0}
\right\|\|w\|^i}{i!} \\
& \leq \bigsum_{i=1}^r \aacc_i \frac{\|w\|^i}{i!}.}
\eeqn
Consider now the possible \sufficient\ termination cases for the 
algorithm and suppose first that termination occurs with 
\accuracy\ as \absolute.  Then, using \req{extDA-difft},
\req{extDA-verif-term-3} and $\omega < 1$, we have that, 
for any $w$ with $\|w\|\leq \delta$,
\vspace*{-2mm}
\beqn{extDA-DTbarDT-bound}
\left|\barDT_r(x,w)- \DT_r(x,w)\right|
\leq \sum_{i=1}^r \aacc_i \frac{\delta^i}{i!}
\leq \omega \xi \frac{\delta^r}{r!}
\leq \xi \frac{\delta^r}{r!}.
\vspace*{-2mm}
\eeqn
If $\barDT_r(x,v) = 0$, we may combine this with \req{extDA-DTbarDT-bound} 
to derive \req{extDA-verif-prop-1}. By contrast, if $\barDT_r(x,v) > 0$,
then since \req{extDA-verif-term-2} failed but \req{extDA-verif-term-3}
holds,
\[
\omega \barDT_r(x,w)
< \sum_{i=1}^r \aacc_i \frac{\delta^i}{i!}
\leq \omega \xi \frac{\delta^r}{r!}.
\]
Combining this inequality with \req{extDA-DTbarDT-bound} yields
\req{extDA-verif-prop-1}. Suppose now that \accuracy\ is \relative.
Then \req{extDA-verif-term-2} holds, and combining it with \req{extDA-difft}
gives that
\vspace*{-2mm}
\[
\left|\barDT_r(x,w)- \DT_r(x,w)\right|
\leq \sum_{i=1}^r \aacc_i \frac{\delta^i}{i!}
\leq \omega \barDT_r(x,w),
\vspace*{-2mm}
\]
for any $w$ with $\|w\|\leq \delta$, which is \req{extDA-verif-prop-2}. 
} 

\noindent
Clearly, the outcome corresponding to our initial aim to obtain a relative
error at most $\omega$ corresponds to the case where \accuracy\ is \relative. 
As we will shortly discover, the two other cases are also needed.

\numsection{The \tal{AR$qp$EDA2} algorithm}\label{ARqpEDA-s}

We now define our adaptive regularization algorithm, which approximately
minimizes the regularized, inexact Taylor-series model defined
\beqn{arqpida-model}
\barm_k(s) \eqdef
\barT_{f,p}(x_k,s) + \frac{\sigma_k}{(p+1)!}\|s\|^{p+1}
\eeqn
at iteration $k$, where $\barT_{f,p}(x,s)$ is described in \req{ipa-barT-def}.
Thus this model uses $p>0$ inexact derivatives, each of which is required
to satisfy bounds of the form \req{eda-vareps-j} for $j\in\ii{p}$. Success or
failure is assessed by comparing the reduction in the Taylor series this gives compared to the 
inexact function value at the resulting trial point. 
Complications arise since optimality
can only be assessed using inexact problem data, and because of
the need for the inexact function values and derivatives to maintain
appropriate coherence with their unknown true values. This makes both the
algorithm and its analysis significantly involved, particularly
since we need to add explicit dynamic-accuracy control to the mix.

Without further ado, and with no further apologies, here is the detailed 
algorithm.
 \newpage

\balgob{extDA-ARqpEDA}{Adaptive-regularization algorithm with 
explicit dynamic \\accuracy  (AR$qp$EDA2)}
{
\vspace*{-2mm}
\begin{description}
\item[Step 0: Initialisation.]
  A criticality order $q$, a degree $p \geq q$, an initial point $x_0\in\Re^n$  
  and an initial regularization parameter $\sigma_0>0$ are given, as well as
  accuracy levels 
  $\epsilon \in (0,1)^q$ and an initial set of absolute derivative accuracies
  $\{\aacc_{i,0}\}_{i=1}^p$.  The constants 
  $\omega$, $\aacc_{\max}$,
  $\varsigma$,
  $\delta_0$, $\theta$, $\eta_1$, $\eta_2$, $\gamma_1$, $\gamma_2$, 
  $\gamma_3$ and $\sigma_{\min}$ are also given and satisfy
  $\aacc_{\max} \geq 0$,
  $\varsigma \in (0,1]$,
  $\theta > 0$, 
  $\delta_0 \in (\epsilon,1]^q$, 
  $\sigma_{\min} \in (0, \sigma_0]$,   \vspace*{-3mm}
  \beqn{arqpeda-eta-gamma}
  0 < \eta_1 \leq \eta_2 < 1, \;\;
  0< \gamma_1 < 1 < \gamma_2 < \gamma_3,
  \eeqn
  \vspace*{-7mm}
  \beqn{extDA-acc-init}
  \omega \in \left(0,\min\left[\half \eta_1,\quarter(1-\eta_2)\right]\right)
  \tim{and}
  \aacc_{i,0} \leq \aacc_{\max} \ms (i\in\ii{p}).
  \eeqn
  Set $k=0$ and $\countacc=0$.  \vspace*{-2mm}

\item[Step 1: Compute the optimality measures and check for termination. ]
  Set $\delta_k^{(0)} = \delta_k$ and evaluate any unavailable
   $\{\overline{\nabla_x^i f}(x_k)\}_{i=1}^p$ 
   to satisfy 
  \beqn{arqpeda-acc-ok}
  \|\overline{\nabla_x^if}(x_k)-\nabla_x^if(x_k)\| \leq \aacc_{i,\countacc}
   \tim{for} i \in \ii{p}.
   \eeqn
   For $j= 1, \ldots ,q$:
  \begin{description}
\item[Step 1.1: ] 
   Compute a displacement $d_{k,j} \in \calB_{\delta_{k,j}}$ such that the
   corresponding Taylor decrement $\barDT_{f,j}(x_k,d_{k,j})$ satisfies
   \beqn{arqpeda-approx-barphi}
   \varsigma \barphi_{f,j}^{\delta_{k,j}}(x_k)\leq \barDT_{f,j}(x_k,d_{k,j}).
   \eeqn
  If the call
  \beqn{arqpeda-verifyo}
  \mbox{\al{CHECK}}\Big(\delta_{k,j},\barDT_{f,j}(x_k,d_{k,j}),
          \{\aacc_{i,\countacc}\}_{i=1}^j, \half \epsilon_j\Big).
  \eeqn
  returns \insufficient, go to Step 5. 
\item[Step 1.2: ] If
  \beqn{arqpeda-termj}
  \barDT_{f,j}(x_k,d_{k,j})
  \leq \frac{\varsigma \epsilon_j}{(1+\omega)} \frac{\delta_{k,j}^j}{j!},
  \eeqn
  consider the next value of $j$.
\item[Step 1.3: ] Otherwise, if
  \beqn{arqpeda-Dmbig}
  \Delta \barm_k(d_{k,j})
  \geq
  \frac{\varsigma\epsilon_j}{2(1+\omega)}\frac{\delta_{k,j}^j}{j!},
  \vspace*{-2mm}
  \eeqn
  go to Step~2 with the index $j_k = j$, and 
  $\delta_k^{(1)} = \delta_k$ and $d_k = d_{k,j}$.
\item[Step 1.4: ]   Otherwise set $\delta_{k,j} = \half \delta_{k,j}$  
  and return to Step~1.1. 
  \end{description}
   Terminate with $x_\epsilon = x_k$  and $\delta_\epsilon = \delta_k$. 
\end{description}
}
{
\vspace*{-2mm}
\begin{description}
\item[Step 2: Step calculation. ]\mbox{}
  \begin{description}
  \item[Step 2.1: ] Compute a step $s_k$ and
  optimality radii $\dsk \in (0,1]^q$ by approximately minimizing the model
  $\barm_k(s)$ from \req{arqpida-model} in the sense that
  \beqn{extDA-descent}
  \Delta \barm_k(s_k) \geq \Delta \barm_k(d_k),
  \eeqn
  \vspace*{-2mm}
   and either
  \beqn{sgtone}
  \|s_k\| \geq 1
  \eeqn
  \vspace*{-2mm}
   or
   \beqn{extDA-mterm}
     \varsigma \phi_{\barm_k,\ell}^{\dskl}(s_k)
     \leq  \DT_{\barm_k,\ell}(s_k,d_{s_k,\ell}^m)  
  \leq \frac{\varsigma\theta(1-\omega)}{(1+\omega)2} 
    \epsilon_\ell \frac{\dskl^\ell}{\ell!}
  \tim{for} \ell \in \ii{q}
  \eeqn
  for some radii $\delta_{s_k} \in (0,1]^q$ and displacements
  $d^m_{s_k,\ell} \in \calB_{\delta_{s_k,\ell}}$.
\item[Step 2.2: ] If the call \vspace*{-4mm}
  \[
  \hspace*{-10.5mm}
  \mbox{\al{CHECK}}\Big( \|s_k\|,\barDT_{f,p}(x_k,s_k),
 \{\aacc_{i,\countacc}\}_{i=1}^p,
  \bigfrac{\varsigma\epsilon_{j_k}}{2(1+\omega)}\,
  \bigfrac{p!\,\left(\delta^{(1)}_{k,j_k}\right)^{j_k}}{j_k!\,
    \max\big[\delta^{(1)}_{k,j_k},\|s_k\|\big]^p}\Big)
  \]
  \vspace*{-4mm}
   or $\|s_k\| <1 $ and any of the calls
  \[
  \hspace*{-7mm}
  \mbox{\al{CHECK}}\Big( \dskl, \DT_{\barm_k,\ell}(s_k,d_{s_k,\ell}^m),
        \{3\max_{t\in\iibe{\ell}{p}} \aacc_{t,\countacc}\}_{i=1}^\ell,
        \frac{\varsigma\theta(1-\omega)\epsilon_\ell}{2(1+\omega)^2} \Big)
       \;\; (\ell \in \ii{q})
  \]
  returns \insufficient, then go to Step~5. 
  \end{description}
\item[Step 3: Acceptance of the trial point. ] 
  Compute $\barf(x_k+s_k)$ ensuring that
  \beqn{extDA-Df+-DT}
  |\barf(x_k+s_k)-f(x_k+s_k)| \leq \omega \barDT_{f,p}(x_k,s_k).
  \eeqn
  Also ensure, either by setting  $\barf(x_k)=\barf(x_{k-1}+s_{k-1})$
  or recomputing $\barf(x_k)$,   that 
  \beqn{extDA-Df-DT}
  |\barf(x_k)-f(x_k)| \leq \omega \barDT_{f,p}(x_k,s_k).
  \eeqn
  Then set
  \vspace*{-2mm}
  \beqn{extDA-rhok-def}
  \rho_k = \frac{\barf(x_k) - \barf(x_k+s_k)}{\barDT_{f,p}(x_k,s_k)}.
  \eeqn
  If $\rho_k \geq \eta_1$, then define
  $x_{k+1} = x_k + s_k$ and $\delta_{k+1} = \delta_{s_k}$ if $\|s_k\|<1$ or
  $\delta_{k+1} = \delta_k^{(1)}$ otherwise.  If $\rho_k < \eta_1$, define
  $x_{k+1} = x_k$ and $\delta_{k+1}=\delta_k^{(1)}$.
\item[Step 4: Regularization parameter update. ]
  Set
  \vspace*{-2mm}
  \beqn{extDA-sigupdate}
  \sigma_{k+1} \in \left\{ \begin{array}{ll}
  {}[\max(\sigma_{\min}, \gamma_1\sigma_k), \sigma_k ]  & \tim{if} \rho_k \geq \eta_2, \\
  {}[\sigma_k, \gamma_2 \sigma_k ]          &\tim{if} \rho_k \in [\eta_1,\eta_2),\\
  {}[\gamma_2 \sigma_k, \gamma_3 \sigma_k ] & \tim{if} \rho_k < \eta_1.
  \end{array} \right.
  \eeqn
  Increment $k$ by one and go to Step~1. 
\end{description}
}
{
\vspace*{-2mm}
\begin{description}
\item[Step 5: Improve accuracy. ]
  For $i\in \ii{p}$, set
  \[
  \aacc_{i,\countacc+1} = \gamma_\aacc\aacc_{i,\countacc},
  \]
  increment $k$ and $\countacc$ by one and return to Step~1 with
  $x_{k+1}=x_k$,  $\delta_{k+1} = \delta_k^{(0)}$ and $\sigma_{k+1}=\sigma_k$.
\end{description}
}

\noindent
Note that extensive use is made of the \al{CHECK} 
algorithm we developed above to ensure that derivative
approximations are sufficiently accurate.
Nonetheless, a number of 
comments are in order to clarify and motivate this extensive 
description.
\begin{itemize}
\item Starting with Step~0, notice that that
we postpone the evaluation of the inexact objective function
$\barf(x_0)$ until Step~3 since Steps~1 and 2 do not depend on its value.

\item Next, examining Step~1, we see that the initialization of Step~1 and Step 1.1
aim at computing, while \req{arqpeda-Dmbig} ensures a lower bound on the model 
decrease, which then guarantees both that $x_k$ is not a global model minimizer 
and also that the first call to \al{CHECK} in Step~2 cannot return 
\absolute\ (see Lemma~\ref{arqpeda-step2} below).  Of course, we need to
show that the loop within Step~1 in which $\delta_{k,j}$ is reduced terminates
finitely, and that the resulting value of $\delta_{k,j}$ is not unduly small (as per
Lemma~\ref{arqpeda-step1loop} below); because of condition
\req{arqpeda-Dmbig}, this entails showing that the model decrease at $d_{k,j}$,
the optimality displacement associated with the violated optimality condition
\req{arqpeda-termj}, is large enough.
We also note that, unless termination occurs, Step~1 specifies the first index
$j\in\ii{q}$ for which $j$-th-order approximate-criticality test
\req{arqpeda-termj} fails.  It is then helpful to distinguish the 
vectors of radii
\beqn{deltass}
\delta_k^{(0)} = \delta_k \tim{at the start of Step 1,}
\eeqn
and inherited from iteration $k-1$, from 
\beqn{deltaes}
\delta_k^{(1)} = \delta_k \tim{at the end of Step 1,}
\eeqn
after possible reductions within that step.  
Clearly, component-wise $\delta_k^{(1)} \leq \delta_k^{(0)}$.

\item Complications arise in Step~2.2 where the step itself and the optimality
displacements associated with model's approximate optimality have to be 
checked for sufficient accuracy.  Just as was the case in Step~1, this entails
possible accuracy improvements, re-evaluation of the derivatives and the need
to recompute the step for the improved model. The absolute accuracy
thresholds passed as the last argument to the two calls to \al{CHECK} 
undoubtedly appear rather mysterious at this stage: the first is 
designed to ensure that an \absolute\ return from \al{CHECK} is impossible,
and the second to ensure an optimality level for the exact problem that is
comparable to that revealed in Step~1. More details will obviously be given in
due course.

\item Step~3 ensures coherence between accuracy on the function
values and accuracy of the model. We stress
that the requirements \req{extDA-Df+-DT} and \req{extDA-Df-DT} do not imply 
that we need to know the true $f$, we only need some mechanism to ensure 
that $x_k$ and $x_k+s_k$ satisfy the required bounds and that are needed to 
guarantee convergence. Again, a new value of $\barf(x_k)$ has to 
be computed to ensure \req{extDA-Df-DT} in Step~3 only when $k > 0$ and
$\barDT_{f,p}(x_{k-1},s_{k-1})>\barDT_{f,p}(x_k,s_k)$, in which case the 
(inexact) function value is computed twice rather than once during that 
iteration. 
As is standard, iteration $k$ is said to be successful when $\rho_k \geq \eta_1$ 
and $x_{k+1}= x_k+s_k$, and we define $\calS$, $\calU$, $\calA$, the sets of
\emph{successful}, \emph{unsuccessful} and \emph{accuracy-improving}
iterations, respectively, as well as $\calT$, $\calS_k$, $\calU_k$, $\calA_k$
and $\calT_k$ by
\beqn{trqida-SU-def}
\begin{array}{rl}
\calS & \eqdef
 \{ k \in \mathbf{N} \mid \mbox{Step~5 is not executed and } 
 \rho_k \geq \eta_1 \},\\
\calU & \eqdef
 \{ k \in \mathbf{N} \mid \mbox{Step~5 is not executed and } 
 \rho_k < \eta_1 \} \tim{and}\\
\calA & \eqdef
 \{ k \in \mathbf{N} \mid \mbox{Step~5 is  executed} \}.
\end{array}
\eeqn
If $\calT \eqdef \calS \cup \calU$, we also define
\beqn{trqida-SUTk-def}
\begin{array}{ll}
\calS_k \eqdef \calS \cap \iiz{k},
& \hspace*{8mm}
\calU_k \eqdef \calU \cap \iiz{k},\\
\calA_k \eqdef \calA \cap \iiz{k}
& \mbox{and} \;\;
\calT_k \eqdef \calT_k \cap \iiz{k},
\end{array}
\eeqn
the corresponding sets up to iteration $k$. Notice that $x_{k+1} = x_k+s_k$
for $k \in \calS$, while $x_{k+1}=x_k$ for $k\in \calU \cup \calA$.
Note also that the objective function is evaluated at most twice per successful
or unsuccessful iteration (i.e.\ once for every $k\in\calT$), and derivatives
are evaluated once per successful or accuracy improving iteration (i.e.\ once
for every $k\in\calS\cup\calA$).

\item Step~4 is the standard regularization parameter update.

\item Finally, Step~5 describes the accuracy improvement mechanism.
\end{itemize}

\noindent
Given the definitions \req{trqida-SUTk-def}, we are now able to show that the
number of iterations in $\calT_k$ is at most a multiple of that in $\calS_k$.

\llem{arqpeda-SvsU}
{
Suppose that the \al{AR$qp$EDA2} algorithm is used and that $\sigma_k \leq
\sigma_{\max}$ for some $\sigma_{\max} > 0$.  Then
\beqn{arqpida-unsucc-neg}
|\calT_k| \leq |\calS_k| \left(1+\frac{|\log\gamma_1|}{\log\gamma_2}\right)+
\frac{1}{\log\gamma_2}\log\left(\frac{\sigma_{\max}}{\sigma_0}\right).
\eeqn
}

\proof{
Observe that $\sigma_{k+1} = \sigma_k$ for $k\in\calA$.
The regularization parameter update \req{extDA-sigupdate}
now gives that, for each $k$,
\[
\gamma_1\sigma_j \leq \max[\gamma_1\sigma_j,\sigma_{\min}]
 \leq \sigma_{j+1}, \ms j \in \calS_k,
\tim{ and }
\gamma_2\sigma_j \leq \sigma_{j+1}, \ms j \in \calU_k.
\]
Thus we deduce inductively that
\[
\sigma_0\gamma_1^{|\calS_k|}\gamma_2^{|\calU_k|}\leq \sigma_k.
\]
Therefore, using our assumption that $\sigma_k\leq \sigma_{\max}$, we
deduce that
\[
|\calS_k|\log{\gamma_1}+|\calU_k|\log{\gamma_2}\leq
\log\left(\frac{\sigma_{\max}}{\sigma_0}\right),
\]
which then implies that
\[
|\calU_k|
\leq - |\calS_k|\frac{\log\gamma_1}{\log\gamma_2}
      + \frac{1}{\log\gamma_2}\log\left(\frac{\sigma_{\max}}{\sigma_0}\right),
\]
since $\gamma_2>1$. The desired result \req{arqpida-unsucc-neg} then follows
from the equality $|\calT_k| = |\calS_k| + |\calU_k|$
and the inequality $\gamma_1 < 1$ given by \req{arqpeda-eta-gamma}.
}

\numsection{Evaluation complexity for the \tal{AR$qp$EDA2} algorithm}

Our analysis of the \al{AR$qp$EDA2} algorithm will be carried out under the
following assumptions.

\begin{description}
\item[f.D0qL: ] $f$ is $q$ times continuously differentiable in $\Re^n$ and,
  for all $j \in \iiz{q}$, the $j$-th derivative of $f$ is Lipschitz continuous
  with Lipschitz constant $L_{f,j}$, that is there exist constants $L_{f,j}
  \geq 1$ such that
  \[
  \|\nabla_x^j f(x) - \nabla_x^j f(y) \| \leq L_{f,j} \|x - y \|
  \]
  for all $x, y \in \Re^n$ and all $j \in \iiz{q}$.

\item[f.Bb: ]There exists a constant $f_{\rm low}$ such that
  $
  f(x) \geq f_{\rm low} \tim{for all } x\in\Re^n.
  $
\end{description}
We then define
\beqn{arqpida-Lf-def}
L_f \eqdef \max_{j\in\iiz{p}} L_{f,j} \geq 1.
\eeqn
We recall that the derivatives of the objective function at the iterates $x_k$ 
remain bounded under \textbf{f.D0pL}. Moreover, the absolute accuracies
$\{\aacc_i\}_{i=1}^p$ never increase in the course of the \al{AR$qp$EDA2}
algorithm, and are initialised so that $\aacc_i\leq \aacc_{\max}$ for all
$i\in\ii{p}$.  As a consequence, \textbf{f.D0pL}, standard error bounds for
Lipschitz functions (see \cite[Corollary~A.8.4]{CartGoulToin21}, for instance),
\req{arqpida-Lf-def} and \req{arqpeda-acc-ok}  imply that, for each  $\ell\in \ii{p}$,
\beqn{arqpeda-dermk-bounded}
\|\overline{\nabla_x^\ell f}(x_k)\|
 \leq \|\nabla_x^\ell f(x_k)\| + 
   \|\overline{\nabla_x^\ell f}(x_k)-\nabla_x^\ell f(x_k)\|
 \leq L_f + \aacc_{\max}
 \eqdef \barL_f.
\eeqn

\subsection{The outcomes of Step~1}

We start by considering the result of performing Step~1 and show that the loop
reducing $\delta_{k,j}$ generated by the possible return to Step~1.1 from
Step~1.4 is finite.  This is done in two stages: we start by expressing 
a general property of the value of the model compared to the truncated Taylor
series, which we will subsequently apply to obtain the desired conclusion.

\llem{arqpida-DmvsDT}
{Suppose that \textbf{f.D0pL} holds, 
that $\barm_k(s)$ is the inexact model \req{arqpida-model}
corresponding to some approximate derivatives 
values $\{\overline{\nabla_x^\ell f}(x_k)\}_{\ell=1}^p$, and 
that $\barL_f$ is given by \req{arqpeda-dermk-bounded}. 
Given $\alpha > 0$ and
\beqn{alpdeltida}
 \delta \in \left(0, \min\!\left(1, \frac{\alpha}{4 \max[\barL_f,\sigma_k]}
 \right)\right],
\eeqn
suppose that there is a displacement $d \in \calB_{\delta}$ for which 
\beqn{arqpida-lDm1}
\barDT_{f,j}(x_k,d) \geq \alpha \frac{\delta^j}{j!}
\eeqn
for some $j\in\ii{q}$. Then
\beqn{arqpida-DMbig1}
\Delta \barm_k(d) 
\geq \half \alpha \frac{\delta^j}{j!}.
\eeqn
}

\proof{
The proof is built on the sequence of inequalities
\[
\begin{array}{rl}
\Delta \barm_k(d_k)
& = \barDT_{f,j}(x_k,d)
     + \bigsum_{\ell=j+1}^p\frac{1}{\ell!}\overline{\nabla_x^\ell f}(x_k)[d]^\ell
     + \bigfrac{\sigma_k}{(p+1)!}\|d\|^{p+1}\\
& \geq \barDT_{f,j}(x_k,d) - 
 \bigsum_{\ell=j+1}^p\bigfrac{\delta^\ell}{\ell!}\|\overline{\nabla_x^\ell f}(x_k)\|
    - \bigfrac{\sigma_k}{(p+1)!}\delta^{p+1}\\
& \geq \alpha \bigfrac{\delta^j}{j!}
     - \max[\barL_f,\sigma_k] \bigsum_{\ell=j+1}^{p+1}\frac{\delta^\ell}{\ell!} \\
& \geq \alpha \bigfrac{\delta^j}{j!}
     - 2 \max[\barL_f,\sigma_k] \frac{\delta^{j+1}}{j!} \\
& =   \bigfrac{\delta^j}{j!}
       \left(\alpha-2\max[\barL_f,\sigma_k]\delta\right)
\end{array}
\]
that arise from the triangle inequality and because of the assumptions made. 
The required bound \req{arqpida-DMbig1} then follows because of 
\req{alpdeltida}.
}

\noindent
We now show that looping inside Step~1 is impossible.

\llem{arqpeda-step1loop}{
Suppose that \textbf{f.D0pL} holds  
 and that algorithm \al{AR$qp$EDA2} has reached the test
\req{arqpeda-Dmbig} in Step~1.3. Suppose also that
\beqn{arqpeda-delta-bound-S1}
\delta_{k,j} \leq  
\frac{\varsigma\epsilon_j}{4(1+\omega) \max[\barL_f, \sigma_k]}
\eeqn
with $\barL_f$ is given by \req{arqpeda-dermk-bounded}.
Then \req{arqpeda-Dmbig} holds, and thus no return from
Step~1.5 to Step~1.2 is possible.
}

\proof{
That the algorithm has reached the test \req{arqpeda-Dmbig} in
Step~1.3 implies that \req{arqpeda-termj} failed  and thus
\[
\barDT_{f,j}(x_k,d_k)
> \left(\frac{\varsigma\epsilon_j}{1+\omega}\right) \frac{\delta_{k,j}^j}{j!}.
\]

Hence, using Lemma~\ref{arqpida-DmvsDT} with 
and the fact that
\req{arqpeda-delta-bound-S1} implies \req{alpdeltida}
for $\delta = \delta_{k,j}$, we deduce that 
\req{arqpeda-Dmbig} must hold from \req{arqpida-DMbig1}.
}

\noindent
Thus the loop within Step~1 is finite, \req{arqpeda-Dmbig} eventually holds
and we may therefore analyze the possible outputs of Step~1.

\llem{arqpeda-step1}{
  Suppose that \textbf{f.D1pL} holds, and that the \al{AR$qp$EDA2} algorithm
  is applied. Then one of three situations may occur at Step~1 of iteration $k$:
  \begin{enumerate}
  \item[(i) ] the \al{AR$qp$EDA2} algorithm terminates with $x_\epsilon$ an
    $(\epsilon,\delta)$-approximate $q$-th-order minimizer, or
  \item[(ii) ] control is passed to Step~5, or
  \item[(iii) ] the call \req{arqpeda-verifyo} returns \relative\ for some
    $j\in\ii{q}$ for which \req{arqpeda-Dmbig} holds for a
    $\delta_{k,j}^{(1)}$ satisfying
    \beqn{arqpeda-delta-lower-S1}
    \delta_{k,j}^{(0)}  \geq \delta_{k,j}  \geq \delta_{k,j}^{(1)}
    \geq  \min\left[ 
    \frac{\varsigma\epsilon_j}{8(1+\omega)
         \max[\barL_f, \sigma_k]}, \delta_{k,j}^{(0)}   \right].
    \eeqn
    and control is then passed to Step~2 with 
  \beqn{arqpeda-good-phi}
  (1-\omega) \barDT_{f,j}(x_k,d_{k,j})
  \leq \phi_{f,j}^{\delta_{k,}^{(0)}}(x_k)
  \leq \left(\frac{1+\omega}{\varsigma}\right)\barDT_{f,j}(x_k,d_{k,j})
  \eeqn
  being satisfied.
  \end{enumerate}
  Moreover, outcome (ii) is impossible whenever
  \beqn{arqpeda-step1-nostep5}
  \max_{i\in\ii{j}}\aacc_i
  \leq
  \frac{\varsigma\omega}{4} 
  \min\left[\frac{\varsigma\epsilon_j}{8(1+\omega)
 \max[\barL_f,\sigma_k]}, \delta_{k,j}^{(0)}\right]^{j-1}
  \frac{\epsilon_j}{j!}.
  \eeqn 
}

\proof{

  Suppose first that branching to Step~5 does not occur in Step~1.2 for any
  $j\in \ii{q}$ (i.e., outcome (ii) does not occur).  This ensures that the
  call \req{arqpeda-verifyo} always returns either \relative\ or \absolute.
  Consider any $j\in \ii{q}$ and notice that Step~1.1  yields
  \req{extDA-rvareps-j} with $T_r=T_{f,j}$ and $r=j$, 
  so that the assumptions of  Lemma~\ref{eda-verify-l} are
  satisfied. Moreover, because of \req{arqpeda-approx-barphi} and the fact that
  $\barphi_{f,j}^{\delta_k}(x_k)\geq0$ by definition,
  we have that   $\barDT_{f,j}(x_k,d_{k,j})\geq 0$.

  If the call \req{arqpeda-verifyo} returns \absolute, 
  then  Lemma~\ref{eda-verify-l}(ii)  with
  $\xi = \half  \varsigma \epsilon_j$, ensures that
  \beqn{trqida-term-j}
  \barDT_{f,j}(x_k,d_{k,j}) \leq \half \varsigma \epsilon_j \bigfrac{\delta_k^j}{j!}
  \leq \frac{\varsigma \epsilon_j}{1+\omega} \; \bigfrac{\delta_k^j}{j!},
  \eeqn
  using the requirement that $\omega < 1$.   Moreover,
  if $d^*_{k,j} \in \calB_{\delta_k}$ is a global maximizer of 
  $T_{f,j}(x_k,d)$ over all $d \in \calB_{\delta_k}$, 
  we may again invoke \req{extDA-verif-prop-1} 
  with $\xi = \half \varsigma\epsilon_j$ 
  together with the triangle inequality to see that
  \beqn{trqida-jopt}
 \arr{rl}{\phi^{\delta_k}_{f,j}(x_k) 
  & = \DT_{f,j}(x_k,d^*_{k,j}) \\*[2ex]
  & \leq \barDT_{f,j}(x_k,d^*_{k,j}) 
     + | \DT_{f,j}(x_k,d^*_{k,j}) - \barDT_{f,j}(x_k,d^*_{k,j}) | \\
  & \leq  \epsilon_j \; \bigfrac{\delta_k^j}{j!}.
  }
  \eeqn
  By contrast, if the call \req{arqpeda-verifyo} returns
  \relative, observe that
  \beqn{trqeda-cor-ineq}
  \varsigma \barDT_{f,j}(x_k,d)
  \leq \varsigma \barphi_{f,j}^{\delta_k}(x_k)
  \leq \barDT_{f,j}(x_k,d_{k,j})
  \eeqn
  for any $d\in \calB_{\delta_k}$ because of \req{arqpeda-approx-barphi},
  and thus we have that
  \[
  \begin{array}{rl}
  \varsigma \DT_{f,j}(x_k,d) 
  & \leq \varsigma \left[\barDT_{f,j}(x_k,d) + 
    \left|\barDT_{f,j}(x_k,d)- \DT_{f,j}(x_k,d)\right|\right] \\*[1ex]
  & \leq \varsigma (1+\omega)\barDT_{f,j}(x_k,d)\\*[1ex]
  & \leq (1+\omega) \barDT_{f,j}(x_k,d_{k,j})
  \end{array}
  \]
  using Lemma~\ref{eda-verify-l}(iii), 
  and the second inequality in \req{arqpeda-good-phi} follows
  by picking $d = d^*_{k,j}$. Similarly
  \[
  \DT_{f,j}(x_k,d) 
  \geq \barDT_{f,j}(x_k,d) - 
   \left|\barDT_{f,j}(x_k,d)- \DT_{f,j}(x_k,d)\right|
  \geq (1-\omega)\barDT_{f,j}(x_k,d)
  \]
  for any $d\in \calB_{\delta_k}$, again using Lemma~\ref{eda-verify-l}(iii).
  Hence
  \[
  \bigmax_{\|d\|\leq \delta_k}\DT_{f,j}(x_k,d) 
   \geq (1-\omega) \bigmax_{\|d\|\leq \delta_k} \barDT_{f,j}(x_k,d)
  \geq (1-\omega)\barDT_{f,j}(x_k,d_{k,j}),
  \]
  which is the first inequality in \req{arqpeda-good-phi}.
  Thus we obtain that, for any $j\in\ii{q}$, either
  both \req{trqida-term-j} and \req{trqida-jopt} hold,
  or \req{arqpeda-good-phi} holds. 

  If the $j$ loop continues to termination, then
  for each   $j\in \ii{q}$, we must have that either 
  \req{trqida-jopt} holds or
  \[
  \phi_{f,j}^{\delta_k}(x_k)
  \leq \left(\frac{1+\omega}{\varsigma}\right) \barDT_{f,j}(x_k,d_{k,j})
  \leq \epsilon_j \frac{\delta_k^j}{j!},
  \]
  where we used \req{arqpeda-good-phi} and the fact that \req{arqpeda-Dmbig}
  must be violated for the loop to continue. Thus 
  in either event \req{trqida-jopt} holds for   all $j \in \ii{q}$, 
  and thus, as \req{optex} holds with $x_\epsilon = x_k$
  and $\delta_\epsilon = \delta_k$, outcome (i) will occur.

  If, instead, control passes to Step~2, we must now show that the conclusions
  of outcome (iii) hold. Firstly, observe that the mechanism of Step~1 ensures
  that the first inequality in \req{arqpeda-delta-lower-S1} is satisfied.
Moreover, if the loop on $j$ does
not finish and branching to Step~5 does not happen, there
must be a $j$ such that \req{arqpeda-termj} is violated and the test
\req{arqpeda-Dmbig} is reached. But we have shown in
Lemma~\ref{arqpeda-step1loop} that the number of times that
\req{arqpeda-Dmbig} is violated---and hence $\delta_{k,j}$ is halved---is
finite since \req{arqpeda-Dmbig} must hold as soon as
\req{arqpeda-delta-bound-S1} holds. At this stage, if one or more
halvings happened, the resulting value of $\delta_{k,j}$ cannot be smaller
than half of the right-hand side of \req{arqpeda-delta-bound-S1}, and so
\req{arqpeda-delta-lower-S1} also holds. If no halving of $\delta_{k,j}$
occurred, $\delta_{k,j}^{(1)} = \delta_{k,j}^{(0)}$ and
\req{arqpeda-delta-lower-S1} obviously holds.
In addition, since \req{arqpeda-Dmbig} ultimately holds for
some $j_k\in \ii{q}$, \req{arqpeda-verifyo} cannot  have 
returned \insufficient\ (by assumption) or \absolute\ (because, in view of
\req{arqpida-model}, \req{trqida-term-j} would have prevented
\req{arqpeda-Dmbig}), and thus it returned \relative. Hence outcome (iii)
occurs.

Finally, Lemma~\ref{eda-verify-l} (i) shows that outcome (ii) cannot happen if
\req{arqpeda-step1-nostep5} holds with  $\delta_k$  being the smallest
$\delta_{k,j}$ that can occur during the execution of Step~1, which is
$\delta_{k,j}^{(1)}$.  In other words, outcome (ii) cannot happen if 
\[
\max_{i\in\ii{j}} \aacc_i
\leq
\frac{\varsigma\omega}{4} \epsilon_j 
\frac{\left(\delta_{k,j}^{(1)}\right)^{j-1}}{j!}.
\]
Substituting \req{arqpeda-delta-lower-S1} into this bound then reveals
\req{arqpeda-step1-nostep5}.
} 

\subsection{The outcomes of Step~2}

Our next task is to investigate Step~2 in more detail. Our aim is 
to show that it is possible to compute a step $s_k$ that satisfies 
\req{extDA-descent} and either \req{sgtone} or \req{extDA-mterm}.
We start our analysis by stating a suitable bound on model decrease, 
which now involves the inexact Taylor series $\barT_{f,p}(x_k,s_k)$
of degree $p$. 

\llem{arqpeda-Dm}{%
The mechanism of the \al{AR$qp$EDA2} algorithm guarantees that
\beqn{extDA-Dphi}
\barDT_{f,p}(x_k,s_k) \geq \frac{\sigma_k}{(p+1)!} \|s_k\|^{p+1}
\eeqn
for all $k \in \calT$,
and \req{extDA-rhok-def} is well-defined.
}

\proof{
Since $k \in \calT$, both \req{arqpeda-Dmbig}
and \req{extDA-descent} must hold at iteration $k$ for some $j \in \ii{q}$.
We then have that
\[
0 < \frac{\varsigma\epsilon_j}{4}\frac{\delta_{k,j}^j}{j!} \leq  
 \Delta \barm_k(d_{k,j})  \leq 
 \Delta \barm_k(s_k)  = \barDT_{f,p}(x_k,s_k)
  -\frac{\sigma_k}{(p+1)!} \|s_k\|^{p+1},
\]
using \req{extDA-descent} which ensures both that $s_k \neq0$ and
\req{extDA-Dphi} holds.
} 

Because the second set of calls to \al{CHECK} in Step~2.2 aim to
check  the accuracy of the Taylor expansion of the \emph{model}, we need to
consider $\{\nabla_d^j T{\!}_{\barm_k,j}(s_k,0)\}_{j=1}^p$ rather than the
$\{\nabla_x^j f(x_k)\}_{j=1}^p$ that we have used so far. It is easy to
verify that these (approximate) derivatives are given by 
\beqn{extDA-ders-mk-plus}
\nabla_d^jT{\!}_{\barm_k,j}(s_k,0)
=\sum_{\ell=j}^p \frac{\overline{\nabla_x^\ell f}(x_k)[s_k]^{\ell-j}}{(\ell-j)!}
+ \frac{\sigma_k}{(p+1)!} \left[\nabla_s^j\|s\|^{p+1}\right]_{s=s_k},
\eeqn
where the last term of the right-hand side is computed exactly.  
This yields the following error bound.

\llem{extDA-abserr-derm-l}{
Suppose that $\|s_k\| \leq 1$.  Then, for all $j \in \ii{p}$,
\beqn{extDA-abserr-derm}
\begin{array}{lcl}
\left|\nabla_d^jT{\!}_{\barm_k,j}(s_k,0)-\nabla_d^jT_{m_k,j}(s_k,0)\right|
& \leq & 3 \,\bigmax_{\ell\in\iibe{j}{p}} \|\overline{\nabla_x^\ell f}(x_k) - \nabla_x^\ell f(x_k) \|\\
& \leq & 3 \,\bigmax_{\ell\in\iibe{j}{p}} \aacc_{\ell,\countacc}.
\end{array}
\eeqn
}

\proof{
Using the triangle inequality, \req{extDA-ders-mk-plus}, the inequality
$\|s_k\| \leq 1$, we have that
\[
\begin{array}{lcl}
\left|\nabla_d^j T_{\barm_k,j}(s_k,0)-\nabla_d^jT_{m_k,j}(s_k,0)\right|
& \!\! \leq \!\! & \bigsum_{\ell=j}^p
         \left|\overline{\nabla_x^jf}(x_k)-\nabla_x^jf(x_k)\right|
         \bigfrac{\|s_k\|^{\ell-j}}{(\ell-j)!}\\
& \!\! \leq \!\! & \bigmax_{\ell\in\iibe{j}{p}}\aacc_{\ell,\countacc} 
 \bigsum_{\ell=j}^p\bigfrac{1}{(\ell-j)!}
\end{array}
\]
for all $j\in\ii{p}$,
and \req{extDA-abserr-derm} then follows from the fact that
\[
\bigsum_{\ell=j}^p\bigfrac{1}{(\ell-j)!}
\leq 1 + \bigsum_{\ell=1}^{p-j}\bigfrac{1}{\ell!}
\leq 1 + e.
\]
}

\noindent
As a consequence, to be safe, we must require three times more accuracy for the
derivatives of the model at $s_k$ than what would be required at $s=0$.

Next, we provide an upper bound on the norm of the step.

\llem{arqpeda-upper-bound-sk}{
Suppose that \textbf{f.D0pL} holds,
that $\barL_f$ is given by \req{arqpeda-dermk-bounded},
and that a step $s_k$ has been found such that \req{extDA-descent} holds.
Then we have that
\beqn{arqpeda-arqp-sk-upper-bd}
\|s_k\|\leq 
\max\left[ \frac{2\barL_f(p+1)!}{\sigma_{\min}}, 
 \left(\frac{2\barL_f(p+1)!}{\sigma_{\min}}\right)^{\sfrac{1}{p}}\right]
\eqdef \kappa_s.
\eeqn
}

\proof{
Using \req{extDA-sigupdate}, \req{extDA-Dphi}, 
the definition \req{ipa-taylor-dec}, the Cauchy-Schwarz inequality, 
and the bounds \req{arqpeda-dermk-bounded} and $\sum_{i=1}^p1/i! < e$, we have that
\[
\arr{rl}{\bigfrac{\sigma_{\min}}{(p+1)!} \|s_k\|^{p+1} 
 & \leq \bigfrac{\sigma_k}{(p+1)!} \|s_k\|^{p+1}
  < \barDT_{f,p}(x_k,s_k) \\
 &
 \leq \bigsum_{i=1}^p \frac{1}{i!}\|\overline{\nabla_x^i f}(x)\|\|s\|^i
   \leq \barL_f \max[\|s_k\|,\|s_k\|^p] \bigsum_{i=1}^p \frac{1}{i!} \\
 &  \leq 2 \barL_f \max[\|s_k\|,\|s_k\|^p],}
\]
which then leads directly to \req{arqpeda-arqp-sk-upper-bd}.
}  

\noindent
We are now in position to elucidate the possible outcomes of Step~2.

\llem{arqpeda-step2}
{Suppose that \textbf{f.D0pL} holds. Then, Step~2 of the
  \al{AR$qp$EDA2} algorithm is well-defined, and either branches to Step~5 
  or produces a pair $(s_k,\dsk)$
  such that
  \beqn{arqpeda-relacc-sk-ok}
  \left|\barDT_{f,p}(x_k,s_k)-\DT_{f,p}(x_k,s_k)\right|
  \leq \omega \barDT_{f,p}(x_k,s_k),
  \eeqn
  and, either
  \begin{itemize}
  \item[(i)]  $\|s_k\| \geq 1$,  or
  \item[(ii)]
  $\|s_k\|<1$ and 
  \beqn{extDA-good-phim}
    \phi_{m_k,j}^{\dskj}(s_k)\leq
    \frac{\theta(1-\omega)}{1+\omega} \,\epsilon_j \,\frac{\dskj^j}{j!}
  \tim{ for } j\in \ii{q}
  \eeqn
  \end{itemize}
  for some $\dsk$ for which 
  \beqn{arqpeda-skeqone}
  \dskj=1 \tim{for} j\in\ii{\min(2,q)}
  \eeqn
  and
  \beqn{arqpeda-kappad-def}
  \dskj\geq \min\left[\kappa_\delta(\sigma_k)\epsilon_j, \delta_{k,j}^{(0)}\right] 
  \tim{for} j\in\iibe{3}{q},
\eeqn
  where
  \beqn{arqpeda-kappads-def}
  \kappa_\delta(\sigma) \eqdef
    \frac{\varsigma\theta(1-\omega)}{8(1+\omega)(3\barL_f+\sigma)} < 1
  \eeqn
  and $\barL_f$ is given by \req{arqpeda-dermk-bounded}.
   Moreover, 
  \beqn{arqpeda-kappad1}
  \delta_{k,j}^{(1)}
  \geq \min\left[\kappa_\delta(\sigma_{k,\max})\,\epsilon_j,\delta_{k,j}^{(0)}\right]
  \eeqn
  for all $k\geq 1$ and $j\in\ii{q}$,
   where $\sigma_{k,\max} = \bigmax_{i\in\iiz{k}}\sigma_i$.
   Finally, $k \in \calT$ whenever
   \beqn{arqpeda-stop-step2}
  \max_{i\in\ii{j}}\aacc_{i,\countacc}
  \leq  \epsmin \min[ \kapsteptwo(\sigma_{k,\max}) \,\epsmin,
  (\delta_{k,\min}^{(0)})]^q
  \eeqn
  where $\kapsteptwo(\sigma)$ is a continuous non-increasing function 
  of $\sigma$, depending only on $\theta$,
  $\varsigma$, $\omega$, $\sigma_{\min}$ and problem constants,
and where
\beqn{arqpeda-delta0min}
\delta_{k,\min}^{(0)}\eqdef \min_{\ell\in\ii{q}}\delta_{k,\ell}^{(0)}.
\eeqn
and $\epsmin\eqdef \min_{j\in\ii{q}} \epsilon_j$.
}

\proof{
The proof proceeds in several stages.

$\bullet$ We first verify that Step~2.1 is well defined. As it turns out, this
conclusion follows from an an almost identical situation arosing in the
analysis of the adaptive regularization algorithm \al{AR$qp$} (using exact
function and derivatives values) when ensuring that a step $s_k$ could be
found for this algorithm. The only significant
difference is that the former this proof used exact derivatives and required
that they were bounded, but now we use approximate ones. Fortunately 
\req{arqpeda-dermk-bounded} provides a substitute bound, and as 
a consequence a straightforward variant of
\cite[Lemma~12.2.9]{CartGoulToin21} still holds 
with $\barL_f$ from \req{arqpeda-dermk-bounded} replacing $L_f$.
In particular, this  lemma (with $\epsilon_j$ replaced by
$\epsilon_j(1-\omega)/(1+\omega)^2$) ensures that \req{extDA-good-phim} is
possible with $\dsk$ satisfying \req{arqpeda-skeqone} and
\req{arqpeda-kappad-def} with \req{arqpeda-kappads-def}. Thus Step~2.1 is well
defined, \req{arqpeda-Dmbig} and \req{extDA-descent} ensure that
\beqn{arqpeda-min-m-decrease}
\Delta \barm_k(s_k)
\geq \Delta \barm_k(d_k)
\geq \frac{\varsigma\epsilon_{j_k}}{2(1+\omega)}\, 
 \frac{\left(\delta_{k,j_k}^{(1)}\right)^{j_k}}{j_k!},
\eeqn
and thus that $s_k\neq 0$, and control passes to~Step 2.2.

$\bullet$ Consider now the evolution of the vector of $\delta_k$ as the
algorithm proceeds. In particular, flow is from Step~1 either to Step~4 (via
Steps~2 and 3) or Step~5 (either directly or via Step~2). The radii are
unaltered on iterations that move via Step 5, and on others, they obey
either \req{arqpeda-delta-lower-S1} or 
\req{arqpeda-skeqone} and \req{arqpeda-kappad-def}.
Therefore, 
\[
\arr{rl}{
  \delta_{k,j}^{(1)} 
  & \geq
  \min \left[
     \bigfrac{\varsigma\theta(1-\omega)\epsilon_j}{8(1+\omega)
  (3\barL_f+\sigma_{k,\max})},
     \bigfrac{\varsigma \epsilon_j}{8(1+\omega)
   \max[\barL_f, \sigma_{k,\max}]}, \delta_{k,j}^{(0)} \right] \\*[2.5ex]
& = \min\left[
    \bigfrac{\varsigma\theta(1-\omega)\epsilon_j}{8(1+\omega)
 (3\barL_f+\sigma_{k,\max})},
    \delta_{k,j}^{(0)}
    \right]
  }
\]
since $\omega \in (0,1)$, which is \req{arqpeda-kappad1}.

$\bullet$ The next step in the proof is to consider the outcome of the 
first call to \al{CHECK} in Step~2.2, and we claim that it can only return 
\insufficient\ or \relative.
Indeed, suppose it returns \absolute. Lemma~\ref{eda-verify-l}(ii)
then implies that
\[
\barDT_{f,p}(x_k,s_k)\leq 
\bigfrac{\varsigma \epsilon_{j_k}}{2(1+\omega)}\,
\left(\bigfrac{p!\,\left(\delta_{k,j_k}^{(1)}\right)^{j_k}}{j_k!\,
 \max\big[\delta_{k,j_k}^{(1)},\|s_k\|\big]^p}\right) 
\frac{\|s_k\|^p}{p!}
\leq
\bigfrac{\varsigma\epsilon_{j_k}}{2(1+\omega)}
\frac{\left(\delta_{k,j_k}^{(1)}\right)^{j_k}}{j_k!}
\]
But the definition of $\barm_k$ and the fact that $s_k\neq 0$ give that
$\barDT_{f,p}(x_k,s_k) > \Delta \barm_k(s_k)$ and hence
\[
\Delta \barm_k(s_k)
< \bigfrac{\varsigma\epsilon_{j_k}}{2(1+\omega)}
 \frac{\left(\delta_{k,j_k}^{(1)}\right)^{j_k}}{j_k!}
\]
which contradicts \req{arqpeda-min-m-decrease}. Our assumption that \al{CHECK}
returned \absolute\ is thus impossible. As a consequence, and unless branching
to Step~5 occurs because it has returned \insufficient, the first call to 
\al{CHECK} in Step~2.2 must return
\relative, which then ensures \req{arqpeda-relacc-sk-ok}
because of Lemma~\ref{eda-verify-l}(iii)

$\bullet$ We now determine a value of the absolute accuracy below 
which the first call to \al{CHECK} in Step~2.2 cannot return
\insufficient.
Recall that the absolute accuracies $\{\aacc_i\}_{i=1}^p$ are reduced by a factor
$\gamma_\aacc$ in Step~5 each time Step~2 passes there.
But \req{extDA-verif-term-2} implies that this will be impossible 
as soon as the accuracies are small enough that
\beqn{arqpeda-nostep5}
\bigsum_{i=1}^p\aacc_{i.\countacc}\bigfrac{\|s_k\|^i}{i!}
\leq
\omega \barDT_{f,p}(x_k,s_k).
\eeqn
But
\[
\barDT_{f,p}(x_k,s_k)
> \Delta \barm_k(s_k)
\geq \Delta \barm_k(d_k) = \Delta \barm_k(d_{k_,j_k})
\geq \frac{\varsigma\epsilon_{j_k}}{2(1+\omega)}\,
 \frac{(\delta_{k,j_k}^{(1)})^{j_k}}{j_k!}
\]
where we successively used the definition of $\barm_k$ \req{arqpida-model},
\req{extDA-descent} and \req{arqpeda-Dmbig}.
Moreover,  the inequality $\delta_{k,j}\leq 1$ and \req{arqpeda-kappad1} 
imply that
\[
\frac{\varsigma\epsilon_j}{2(1+\omega)}\,\frac{(\delta_{k,j}^{(1)})^j}{j!}
\geq \frac{\varsigma\epsilon_j}{2(1+\omega)}\,\frac{(\delta_{k,j}^{(1)})^q}{q!}
\geq \frac{\varsigma\epsilon_j}{2(1+\omega)}\,
    \frac{\min[\kappa_\delta(\sigma_{k,\max})\epsilon_j,\delta_{k,j}^{(0)}]^q}{q!}
\]
so that
\beqn{arqpeda-nostep5-2}
\barDT_{f,p}(x_k,s_k)
\geq\frac{\varsigma\epsilon_{j_k}}{2(1+\omega)}\,
    \frac{\min[\kappa_\delta(\sigma_{k,\max})\epsilon_{j_k},\delta_{k,j_k}^{(0)}]^q}{q!}.
\eeqn
In addition, we know from Lemma~\ref{arqpeda-upper-bound-sk} and the
sum bound $\sum_{i=1}^p1/i!< e$ that
\[
\bigsum_{i=1}^p\aacc_{i.\countacc}\bigfrac{\|s_k\|^i}{i!}
\leq 2 \max[ 1, \kappa_s^p] \max_{i\in\ii{p}}\aacc_{i.\countacc},
\]
Combining this inequality with \req{arqpeda-nostep5-2}, we deduce that
\req{arqpeda-nostep5} holds, and consequently the first call to \al{CHECK} in
Step~2.2 returns \relative\ (as we have just verified it cannot return 
\absolute), as soon as 
\beqn{arqpeda-bV1}
\max_{i\in\ii{p}}\aacc_{i.\countacc}
\leq \omega\frac{\varsigma\epsilon_{j_k}}{4(1+\omega)}\,
     \frac{1}{ \max[ 1, \kappa_s^p]}
     \frac{\min\big[\kappa_\delta(\sigma_{k,\max})\epsilon_{j_k},
 \delta_{k,j_k}^{(0)}\big]^q}{q!}.
\eeqn
Note that this conclusion is independent of $\|s_k\|$.

$\bullet$ We next consider what can be said of the optimality conditions on the
model when $\|s_k\|<1$, and start by finding values of the absolute accuracy
that are acceptably small to prevent any of the second set of calls to
\al{CHECK} in Step~2.2 from returning \insufficient. 
Lemmas~\ref{eda-verify-l}(i) and \ref{extDA-abserr-derm-l}---notice
that the latter ensures that \req{extDA-rvareps-j} holds when we invoke 
the former---and the argument list for \al{CHECK} 
suggest that this is impossible if
\beqn{verify2-sufficient}
3\max_{t\in\iibe{\ell}{p}} \aacc_{t,\countacc} \sum_{i=1}^{\ell} 
\frac{\dskl^i}{i!} \leq
\varsigma\omega \frac{\theta(1-\omega)}{2(1+\omega)^2}\,\epsilon_\ell 
 \frac{\dskl^\ell}{\ell!}
\eeqn
for all $\ell \in \ii{q}$.
But since $\dsk \leq 1$, the sum bound $\sum_{i=1}^p1/i!< e$ again implies that 
\beqn{arqpedq-3-to-6}
3\max_{t\in\iibe{\ell}{p}} \aacc_{t,\countacc} \sum_{i=1}^{\ell} \frac{\dskl^i}{i!} 
\leq 6\max_{t\in\ii{p}} \aacc_{t,\countacc},
\eeqn
while \req{arqpeda-skeqone} and \req{arqpeda-kappad-def}, the
requirement that $\epsilon_j \leq 1$ 
and the bounds $\kappa_\delta(\sigma_{k,\max}) < 1$ and $\ell \leq q$ 
give that
\[
\begin{array}{rl}
\varsigma \omega \bigfrac{\theta(1-\omega)}{2(1+\omega)^2}\,\epsilon_\ell 
\bigfrac{\dskl^\ell}{\ell!} 
& \geq \varsigma \omega \bigfrac{\theta(1-\omega)\epsilon_\ell}{2(1+\omega)^2}\,
 \bigfrac{\min[\kappa_\delta(\sigma_k)\epsilon_\ell,\delta_{k,\ell}^{(0)}]^\ell}{\ell!}\\
& \geq \varsigma \omega \bigfrac{\theta(1-\omega)\epsilon_\ell}{2(1+\omega)^2}\,
 \bigfrac{\min[\kappa_\delta(\sigma_{k,\max})\epsilon_\ell,\delta_{k,\ell}^{(0)}]^\ell}{\ell!}\\
\end{array}
\]
Therefore, in view of \req{arqpedq-3-to-6},
\req{verify2-sufficient} holds, and
\al{CHECK} returns either \absolute\ or \relative, whenever
\beqn{arqpeda-bV2}
\max_{t\in\ii{p}} \aacc_{t,\countacc} \leq 
\varsigma \omega \bigfrac{\theta(1-\omega)\epsmin}{12(1+\omega)^2}\,
 \bigfrac{\min\big[\kappa_\delta(\sigma_{k,\max})\epsmin,\delta_{k,\min}^{(0)}\big]^q}{q!}
\eeqn
where $\delta_{k,\min}^{(0)}$ is given by \req{arqpeda-delta0min}.

Consider the $\ell$-th such call, and note that
\beqn{arqpeda-split}
\DT_{m_k,\ell}(s_k,d)
\leq \DT_{\barm_k,\ell}(s_k,d) + |\DT_{m_k,\ell}(s_k,d)-\DT_{\barm_k,\ell}(s_k,d)|
\eeqn
for any $d \in \calB_{\delta_{s_k,\ell}}$
because of the triangle inequality. If \al{CHECK} returns \absolute, then
\req{arqpeda-split}, Lemma~\ref{eda-verify-l} (ii) 
and the fact that $\varsigma \leq 1$ yield that 
\[
\phi_{m_k,\ell}^{\delta_{s_k,\ell}}(s_k)
  = \bigmax_{\|d\|\leq \delta_{s_k,\ell}}\DT_{m_k,\ell}(s_k,d)
\leq 2\,\varsigma \bigfrac{\theta(1-\omega) \epsilon_\ell}{2(1+\omega)^2}\,
     \bigfrac{\delta_{s_k,\ell}^\ell}{\ell!}
\]
and \req{extDA-good-phim} holds. By contrast, if \al{CHECK} returns
\relative, then \req{arqpeda-split}, Lemma~\ref{eda-verify-l}(iii)
  and \req{extDA-mterm} ensure that, for some $d_{k,\ell}^*\in \calB_{\delta_{s_k,\ell}}$,
\[
\begin{array}{rl}
\varsigma \phi_{m_k,\ell}^{\dskl}(s_k) 
& \leq \varsigma \,(1+\omega)\,\DT_{\barm_k,\ell}(s_k,d_{k,\ell}^*)\\*[1.5ex]
& \leq (1+\omega)\,\varsigma \,\phi_{\barm_k,\ell}^{\dskl}(s_k)\\*[1.5ex]
& \leq \varsigma \bigfrac{\theta(1-\omega)\epsilon_\ell}{(1+\omega)}
      \,\bigfrac{\delta_{s_k,\ell}^\ell}{\ell!}
\end{array}
\]
and \req{extDA-good-phim} holds again.  Thus \req{extDA-good-phim} holds in
both cases.

$\bullet$ We conclude from the above discussion, in particular from
\req{arqpeda-bV1} and \req{arqpeda-bV2}, that Step~2 terminates 
with a pair $(s_k,\dsk)$, for which \req{extDA-good-phim} holds 
if $\|s_k\| \leq 1$, whenever
\[
\max_{i\in\ii{p}}\aacc_{i,\countacc}
\leq \min\big[\kapsteptwo(\sigma_{k,\max})\,\epsmin,\delta_{k,\min}^{(0)}\big]^q\,\epsmin
\]
where
\[
\kapsteptwo(\sigma) \eqdef
\frac{\varsigma\omega \left(\kappa_\delta(\sigma)\right)^q}{4q!(1+\omega)}
\min\left[\bigfrac{1}{\max[1,\kappa_s^p]},
\frac{\theta(1-\omega)}{3(1+\omega)}
\right].
\]
The algorithm then proceeds to Step~3, and thus $k \in \calT$.
} 

\noindent
We bring together two results obtained so far regarding bounds on
$\delta_{k,j}$, namely those in Lemma~\ref{arqpeda-step1} and 
\ref{arqpeda-step2}.

\llem{arqpeda-newkappadelta}{
  Suppose that \textbf{f.D0pL} holds, that iteration $k$ 
  of the \al{AR$qp$EDA2} algorithm is successful, and that
  $\|s_k\|<1$.  Then 
  \[
  \delta_{k+1,j}^{(0)} = 1 \tim{for} j\in\ii{\min[2,q]}
  \]
  and
  \[
  \delta_{k+1,j}^{(0)} \geq \min\left[\kappa_\delta(\sigma_k) \, \epsilon_j,
    \delta_{k,j}^{(0)}\right]
  \tim{for} j \in \iibe{3}{q},
  \]
  where $\kappa_\delta(\sigma)$ is defined in \req{arqpeda-kappads-def}.
}

\proof{
If iteration $k$ is successful, Step~3 sets
$\delta_{k+1,j}^{(0)} = \dsk$. The stated bound then follows from
\req{arqpeda-skeqone} and \req{arqpeda-kappad-def}.
} 

\llem{arqpeda-newallkappadelta}{
  Suppose that \textbf{f.D0pL} holds, and that the \al{AR$qp$EDA2} algorithm 
  does not terminate at (or before) iteration $k$. Then
  \beqn{unco-radiizero}
  \delta_{k+1,j}^{(0)} \geq \kappa_\delta(\sigma_{k,\max}) \, \epsilon_j
  \tim{for} j \in \ii{q},
  \eeqn
  where $\kappa_\delta(\sigma)$ is defined in \req{arqpeda-kappads-def}
  and $\sigma_{k,\max} = \max_{i\in\iiz{k}}\sigma_i$.
  Moreover, $k\in \calT$ whenever
  \beqn{arqpeda-step2-accuracy}
  \min_{i\in\ii{p}}\aacc_{i,\countacc}
  \leq \kapsteptwo (\sigma_{k,\max})\epsmin^{q+1}.
  \eeqn
}

\proof{We prove this by induction over $k$ for each $j \in \ii{q}$.
By assumption, the $j$-th initial radius satisfies
\[
\delta^{(0)}_{0,j} = \delta_{0,j} \geq \epsilon_j
> \kappa_\delta(\sigma_{0}) \, \epsilon_j
= \,\kappa_\delta(\sigma_{0,\max}) \, \epsilon_j.
\]
Suppose that $k = 0$.
Then, as the algorithm does not terminate during this iteration,
Lemma~\ref{arqpeda-step1} indicates that either control is passed
to Step~5, in which case
\beqn{unco-one-to-five}
\delta^{(0)}_{1,j} = 
\delta^{(0)}_{0,j} \geq \kappa_\delta(\sigma_{0,\max})  \, \epsilon_j
\eeqn
as above, or \req{arqpeda-delta-lower-S1} holds, i.e.,
\beqn{unco-deltze}
 \delta_{0,j}^{(1)} \geq  \min\left[ \frac{\varsigma
  \epsilon_j}{8(1+\omega)\max[\barL_f, \sigma_0]},
    \delta_{0,j}^{(0)} \right] \geq \kappa_\delta(\sigma_{0,\max}) \, \epsilon_j.
\eeqn
Step~2 may then pass control to Step~5 (with the same outcome 
\req{unco-one-to-five} as for Step~1), 
but if not Step~3 will either result in
\[\delta_{1,j}^{(0)}
\geq \min\left[\kappa_\delta(\sigma_0)\,\epsilon_j,\delta_{0,j}\right]
\geq \min\left[\kappa_\delta(\sigma_{0,\max})\,\epsilon_j,\delta_{0,j}\right]
\geq \kappa_\delta(\sigma_{0,\max})
 \]
as per Lemma~\ref{arqpeda-newkappadelta} when the 
iteration is successful and $\|s_0\|<1$, or otherwise 
\[
\delta_{1,j}^{(0)} =
 \delta_{0,j}^{(1)} \geq \kappa_\delta(\sigma_{0,\max}) \, \epsilon_j
\]
because of \req{unco-deltze}. Thus in all cases \req{unco-radiizero} 
holds for $k=0$.

Now suppose that \req{unco-radiizero} holds up to iteration $k-1$, i.e.,
\beqn{unco-radiizerok}
 \delta_{k,j}^{(0)} \geq \kappa_\delta(\sigma_{k-1,\max}) \, \epsilon_j
 \geq \kappa_\delta(\sigma_{k,\max}) \, \epsilon_j.
\eeqn
We show it also holds for iteration $k$.
The proof is essentially identical to the $k=0$ case.
Once again, as the algorithm does not terminate during iteration $k$,
Lemma~\ref{arqpeda-step1} shows that either control is passed
to Step~5, in which case 
\beqn{unco-one-to-fivek}
\delta^{(0)}_{k+1,j} = \delta^{(0)}_{k,j} 
\geq \kappa_\delta(\sigma_{k,\max}) \, \epsilon_j
\eeqn
from \req{unco-radiizerok},
or \req{arqpeda-delta-lower-S1} holds, i.e.,
\beqn{unco-deltzek}
\begin{array}{rl}
 \delta_{k,j}^{(1)}  & \geq  
 \min\left[ \bigfrac{\varsigma\epsilon_j}{8(1+\omega)
 \max[\barL_f,\sigma_k]},\delta_{k,j}^{(0)} 
 \right]  \\*[2.5ex]
 & \geq
 \min\left[ \bigfrac{\varsigma}{8(1+\omega)
 \max[\barL_f,\sigma_k]},\kappa_\delta(\sigma_{k,\max}) 
\right] \, \epsilon_j 
 = \kappa_\delta(\sigma_{k,\max}) \epsilon_j
 \end{array}
\eeqn
using \req{unco-radiizerok}. As before, Step~2 may then pass control to Step~5
(and thus \req{unco-one-to-fivek} holds), but if not Step~3 will either result in
\[
\delta_{k+1,j}^{(0)}
\geq \min\left[\kappa_\delta(\sigma_{k}) \, \epsilon_j, \delta_{k,j}^{(0)} \right]
\geq \kappa_\delta(\sigma_{k,\max}) \, \epsilon_j
\]
when the iteration is successful and $\|s_k\|<1$,
using Lemma~\ref{arqpeda-newkappadelta}
and \req{unco-radiizerok}, or otherwise 
\[
\delta_{k+1,j}^{(0)} = \delta_{k,j}^{(1)} \geq \kappa_\delta(\sigma_{k,\max}) 
 \, \epsilon_j
\]
because of \req{unco-deltzek}. This completes the induction
  proving \req{unco-radiizero}. The inequality \req{arqpeda-step2-accuracy}
  then follows by using \req{unco-radiizero} at iteration $k$ (instead of
  $k+1$) in \req{arqpeda-stop-step2}.
}

\subsection{The final complexity bound}

We are now poised to consider the evaluation complexity analysis of the
fully-formed \al{AR$qp$EDA2} algorithm. Having stated the necessary
inexact variant of the standard bound on model decrease in 
\req{extDA-Dphi}, we next show that the regularization parameter $\sigma_k$ 
generated by the algorithm remains bounded, even when $f$ and its derivatives
are computed inexactly. 

\llem{extDA-sigmaupper-lemma}{
Suppose that \textbf{f.D0pL} holds.  Then algorithm \al{AR$qp$EDA2}  ensures that
\beqn{extDA-sigmaupper}
\sigma_k \leq \sigma_{\max} \eqdef
 \max\left[ \sigma_0,\gamma_3 
 \frac{4L_{f,p}}{1-\eta_2}\right]
\vspace*{-2mm}
\eeqn
for all $k\geq 0$.
}

\proof{
Suppose that $k \in \calT$, and that
\beqn{extDA-siglarge}
\sigma_k \geq \frac{4L_{f,p}}{1-\eta_2}
\eeqn
Then the triangle inequality, \req{arqpeda-relacc-sk-ok} that is guaranteed
by Lemma~\ref{arqpeda-step2} as iteration $k$ proceeds to Step~3, 
and \req{extDA-Df-DT} together show that
\[
\begin{array}{lcl}
|\barT_{f,\ell}(x_k,s_k) - T_{f,\ell}f(x_k,s_k)|
& \leq &  |\barf(x_k)-f(x_k)| + |\barDT_{f,\ell}(x_k,s_k)-\DT_{f,\ell}(x_k,s_k)|\\*[1.5ex]
& \leq &  2 \omega |\barDT_{f,\ell}(x_k,s_k)|.
\end{array}
\]
Therefore, again using the triangle inequality, \req{extDA-Df+-DT},
standard error bounds for Lipschitz functions (see \cite[Corollary~A.8.4]{CartGoulToin21}),
\req{extDA-Dphi}, \req{extDA-acc-init} and \req{extDA-siglarge}, we deduce that
\[
\begin{array}{ll}
|\rho_k - 1|
& \leq \bigfrac{|\barf(x_k+s_k) - \barT_{f,p}(x_k,s_k)|}
        {\barDT_{f,p}(x_k,s_k)} \\*[2ex]
& \leq \bigfrac{1}{\barDT_{f,p}(x_k,s_k)}
   \left[|\barf(x_k+s_k) - f(x_k+s_k)| \right.\\*[2ex]
&\hspace*{13.6mm}  \left. + |f(x_k+s_k)-T_{f,p}(x_k,s_k)|
                 + | \barT_{f,p}(x_k,s_k)-T_{f,p}(x_k,s_k)| \right] \\*[2ex]
& \leq \bigfrac{1}{ \barDT_{f,p}(x_k,s_k)}
   \Big[ |f(x_k+s_k)-T_{f,p}(x_k,s_k)|+3\omega \barDT_{f,p}(x_k,s_k)\Big]\\*[2ex]
& \leq \bigfrac{1}{\barDT_{f,p}(x_k,s_k)}
   \left[ \bigfrac{L_{f,p}}{(p+1)!} \|s_k\|^{p+1}\right]
   + 3 \omega \\*[2ex]
& < \bigfrac{L_{f,p}}{\sigma_k} + \bigfrac{3(1-\eta_2)}{4} \\*[2ex]
& \leq 1-\eta_2
\end{array}
\]
and thus that $\rho_k \geq \eta_2$. Then iteration $k$ is very successful in 
that
$\rho_k \geq \eta_2$ and, because of \req{extDA-sigupdate}, 
$\sigma_{k+1}\leq \sigma_k$.  
As a consequence, the mechanism of the algorithm ensures that
\req{extDA-sigmaupper} holds for all $k \in \calT$, while if $k \in \calA$,
Step~5 fixes $\sigma_{k+1} = \sigma_k$.
}

\noindent
We now state a useful technical inequality.

\llem{extDA-phinext}{
  Suppose that \textbf{f.D0pL} holds. Suppose also that iteration $k$ of the
  \al{AR$qp$EDA2} algorithm is successful, that $\|s_k\|<1$, 
  but that iteration $k_+$ is the first iteration after $k$ that
  proceeds to Step~2 rather than Step~5.  Then there exists 
  $j \in \ii{q}$ such that 
  \beqn{extDA-phinext-bound}
  \frac{\varsigma (1-\theta)(1-\omega)}{1+\omega}\,\epsilon_j \,
   \bigfrac{\left(\delta_{k_+,j}^{(0)}\right)^j}{j!}
  \leq  (L_{f,p} + \sigma_{\max})
  \bigsum_{\ell=1}^j 
 \bigfrac{\left(\delta_{k_+,j}^{(0)}\right)^\ell}{\ell!}\|s_k\|^{p-\ell+1},
  \eeqn
where $\sigma_{\max}$ is defined by \req{extDA-sigmaupper}.
}

\proof{
  Observe that our assumption that iteration $k_+$ passes to Step~2 means that
  it does not terminate in Step~1, nor does it pass to Step~5. The latter
  implies that 
  \[
  \barDT_{f,j}(x_{k_+},d_{k_+,j})
  > \left(\frac{\varsigma\epsilon_j}{1+\omega}\right)
    \frac{\left(\delta_{k_+,j}^{(0)}\right)^{j}}{j!}
  \]
  must hold for some $j\in \ii{q}$. Furthermore, since any iteration
  that might occur between $k$ and $k_+$ must have passed through Step~5,
  $\delta_{k_+,j}^{(0)} = \dskj$ as the radii update rule in Step~5 does not 
  change the input value of $\delta_k$. Hence \req{arqpeda-good-phi}
  from Lemma~\ref{arqpeda-step1} ensures that
  \beqn{extDA-no-term}
  \phi_{f,j}^{\delta_{k_+,j}^{(0)}}(x_{k_+})
  >
  \varsigma\left(\frac{1-\omega}{1+\omega}\right) \,\epsilon_{j}\,
  \frac{\left(\delta_{k_+,j}^{(0)}\right)^{j}}{j!}.
  \eeqn
  and therefore that
  \beqn{unco-arqp-longs1}
\begin{array}{l}
  \varsigma\left(\bigfrac{1-\omega}{1+\omega}\right) \,
  \epsilon_j\bigfrac{(\delta_{+,j}^{(0)})^j}{j!} < \phi_{f,j}^{\delta_{+,j}^{(0)}}(x_{k+1})
 =  - \bigsum_{\ell=1}^j\frac{1}{\ell!}\nabla_x^\ell f(x_{k+1})[d]^\ell \\*[3ex]
\hspace*{10mm}
 =  - \bigsum_{\ell=1}^j\frac{1}{\ell!}\nabla_x^\ell f(x_{k+1})[d]^\ell
         + \bigsum_{\ell=1}^j \bigfrac{1}{\ell!}\nabla^\ell_sT_{f,p}(x_k,s_k)[d]^\ell\\*[3ex]
\hspace{20mm} - \bigsum_{\ell=1}^j \bigfrac{1}{\ell!}\nabla^\ell_sT_{f,p}(x_k,s_k)[d]^\ell
- \frac{\sigma_k}{(p+1)!} \bigsum_{\ell=1}^j \bigfrac{1}{\ell!}
 \nabla^\ell_s \left( \|s_k\|^{p+1}\right) [d]^\ell\\*[3ex]
\hspace{20mm} 
+ \bigfrac{\sigma_k}{(p+1)!} \bigsum_{\ell=1}^j \bigfrac{1}{\ell!}
 \nabla^\ell_s \left( \|s_k\|^{p+1}\right) [d]^\ell.
\end{array}
\eeqn
To bound the terms on the right-hand side of \req{unco-arqp-longs1},
firstly, using standard approximation properties for Lipschitz function (once
more, see \cite[Corollary~A.8.4]{CartGoulToin21}),
\beqn{unco-arqp-skd}
\begin{array}{l}
- \bigsum_{\ell=1}^j\frac{1}{\ell!}\nabla_x^\ell f(x_{k+1})[d]^\ell
+ \bigsum_{\ell=1}^j \bigfrac{1}{\ell!}\nabla^\ell_sT_{f,p}(x_k,s_k)[d]^\ell\hspace*{30mm}\\*[2ex]
\hspace{20mm} \leq \bigsum_{\ell=1}^j \bigfrac{\delta_{k+1,j}^\ell}{\ell!}
                   \big\|\nabla_x^\ell f(x_{k+1})-\nabla^\ell_sT_{f,p}(x_k,s_k)\big\|\\*[2ex]
\hspace{20mm} \leq L_{f,p}\bigsum_{\ell=1}^j
                   \bigfrac{\delta_{k+1,j}^\ell}{\ell!(p-\ell+1)!}\|s_k\|^{p-\ell+1},\\*[2ex]
\end{array}
\eeqn
where $L_{f,p}$ is defined in \textbf{f.D1pL}.
Furthermore, because of \req{arqpida-model} and \req{extDA-mterm},
\beqn{unco-arqp-longs6}
\begin{array}{l}
  - \bigsum_{\ell=1}^j \bigfrac{1}{\ell!}\nabla^\ell_sT_{f,p}(x_k,s_k)[d]^\ell
  - \frac{\sigma_k}{(p+1)!} \bigsum_{\ell=1}^j \bigfrac{1}{\ell!}
   \nabla^\ell_s \left( \|s_k\|^{p+1}\right) [d]^\ell\\*[2ex]
\hspace*{15mm} \leq \phi_{m_k,j}^{\dskj}(s_k)
\leq \theta \epsilon_j \,\bigfrac{\dskj^j}{j!}
= \theta \epsilon_j \,\bigfrac{\delta_{k+1,j}^j}{j!},
\end{array}
\eeqn
where the last equality follows as $\delta_{s_k,j}=\delta_{+,j}^{(0)}$ 
if iteration $k$ is successful. Moreover, in view of 
the form of the derivatives of the regularization term $\|s\|^p$ (see
\cite[Lemma~B.4.1]{CartGoulToin21} with $\beta=1$) and
Lemma~\ref{extDA-sigmaupper-lemma}, we also have that
\beqn{unco-arqp-longs5}
\begin{array}{rl}
\bigfrac{\sigma_k}{(p+1)!} \bigsum_{\ell=1}^j \bigfrac{1}{\ell!}
 \nabla^\ell_s \left( \|s_k\|^{p+1}\right) [d]^\ell  
& \leq \bigfrac{\sigma_k}{(p+1)!} \bigsum_{\ell=1}^j \bigfrac{1}{\ell!}
 \left\| \nabla^\ell_s \left( \|s_k\|^{p+1}\right)\right \| \|d\|^\ell \\
& = \sigma_k \bigsum_{\ell=1}^j
 \bigfrac{\|s_k\|^{p-j+1} \|d\|^\ell}{\ell!(p-j+1)!}
\\
& \leq \sigma_{\max} \bigsum_{\ell=1}^j
 \bigfrac{\delta_{k+1,j}^\ell}{\ell!(p-j+1)!}\|s_k\|^{p-j+1}. 
\end{array}
\eeqn
   We then observe that, 
  \beqn{extDA-phim-b}
  \phi_{m_k,j}^{\delta_{k_+,j}^{(0)}}(s_k)
  = \phi_{m_k,j}^{\dskj}(s_k)
  \leq \varsigma \theta\left(\frac{1-\omega}{1+\omega}\right)
     \epsilon_{j}\frac{\dskj^{j}}{j!},
  \eeqn
  where we have used the fact that $\delta_{k_+,j}^{(0)} = \dskj$ for
  $k\in\calS$ to derive the first inequality and the fact that $\|s_k\|<1$  to
  apply the bound \req{extDA-good-phim} of Lemma~\ref{arqpeda-step2}. The
  proof is then concluded by replacing the second inequality in
  \req{unco-arqp-longs6} by \req{extDA-phim-b} and combining the result with
  \req{unco-arqp-longs1}, \req{unco-arqp-skd} and \req{unco-arqp-longs5}.
  
} 

\noindent
Our next step is to provide a lower bound on the step at iterations before
termination.

\llem{arqpeda-longs-l}{
Suppose that \textbf{F.D0pL} holds, iteration $k$ is successful, and that the
\al{AR$qp$EDA2} algorithm does not terminate at first iteration after $k$ that
proceeds to Step~2 rather than Step~5. Suppose also
that the algorithm chooses $\delta_k$ for each $k$ such that the conclusions of
Lemma~\ref{arqpeda-newkappadelta} hold. Then there exists 
a $j\in\ii{q}$ such that
\beqn{extDA-longs}
\|s_k\|\geq
\left\{\begin{array}{ll}
  \left(\bigfrac{\varsigma(1-\theta)(1-\omega)}
    {2j!(L_{f,p}+\sigma_{\max})(1+\omega)}\right)^{\frac{1}{p-j+1}}
    \,\epsilon_j^{\frac{1}{p-j+1}}
   & \tim{if} q \in \{1,2\}, \\*[3.5ex]
    \left(\bigfrac{\varsigma(1-\theta)(1-\omega)\kappa_{\delta,\min}^{j-1}}
    {2j!(L_{f,p}+\sigma_{\max})(1+\omega)}\right)^{\frac{1}{p}}
    \,\epsilon_j^{\frac{j}{p}}
   & \tim{if} q > 2,
\end{array}\right.
\eeqn
where $\kappa_{\delta,\min} \eqdef \kappa_\delta(\sigma_{\max})$ 
and $\sigma_{\max}$ is defined by \req{extDA-sigmaupper}.
}

\proof{
Either $\|s_k\|\geq 1$ (case (i) in Lemma~\ref{arqpeda-step2})  and
\req{extDA-longs} holds automatically because the two lower bounds on its
right-hand side are less than one.  So suppose instead that  $\|s_k\| \leq 1$.
Because the algorithm does not terminate, Lemma~\ref{extDA-phinext}
ensures that \req{extDA-phinext-bound} holds for some $j\in \ii{q}$.  It
is easy to verify that this inequality is equivalent to
\beqn{ineq-ab}
\tau \,\epsilon_j \, (\delta_{+,j}^{(0)})^j
\leq \|s_k\|^{p+1} \chi_j\left(\frac{\delta_{+,j}^{(0)}}{\|s_k\|}\right)
\eeqn
where the function $\chi_j$ is defined by
\beqn{unco-chi-def}
\chi_j(t) = \sum_{\ell=1}^j \frac{t^\ell}{\ell !}
\eeqn
and where we have set
\[
\tau = \frac{(1-\omega)(1-\theta)}{j!(1+\omega)(L_{f,p}+\sigma_{\max})}.
\]
In particular, since $\chi_j(t) \leq 2 t^j$ for $t \geq 1$, we have that
\beqn{ineq-good}
\tau \,\epsilon_j
\leq 2 \|s_k\|^{p+1}\left(\frac{1}{\|s_k\|}\right)^j
= 2 \|s_k\|^{p-j+1}
\eeqn
whenever $\|s_k\| \leq \delta_{k+1,j}$.

If $j$ is 1 or 2, by assumption $\delta_{+,j}^{(0)}=1$ and 
$\|s_k\| \leq 1 = \delta_{+,j}^{(0)}$.  Thus
\req{ineq-good} yields the first case of \req{extDA-longs}.
Otherwise, if $j \geq 3$, by assumption \req{unco-radiizero} holds. 
In this case, if $\|s_k\|\leq \delta_{+,j}^{(0)}$, we may again deduce
from \req{ineq-good} that the first case of \req{extDA-longs} holds, and
this implies that the second case also holds since
$\kappa_{\delta,\min}<1$ and\footnote{
This is easily verified by noting that, if $\lambda(t) = t/(j(t-j+1))$, then
$\lambda(j)=1$ and $\lambda^{(1)}(t)\leq 0$ for all $t \geq j$.
}
$1/(p-j+1)\leq j/p$ for $p\geq j\geq 1$.
Consider therefore the alternative for which $\|s_k\| > \delta_{+,j}^{(0)}$.  Then
\req{ineq-ab}, and noting that $\chi_j(t) < 2t$ for $t\in [0,1]$,
we deduce that
\[
\tau \,\epsilon_j \, (\delta_{+,j}^{(0)})^j
\leq 2 \|s_k\|^{p+1} \left(\frac{\delta_{+,j}^{(0)}}{\|s_k\|}\right),
\]
which, with \req{unco-radiizero}, implies the second case of 
\req{extDA-longs}. 
}

\noindent
We now consolidate the above results by stating lower
bounds on the minimal model and function decreases at successful iterations.

\llem{lem-arqpeda-model-decrease}{
Suppose that \textbf{f.D1qL} holds. Then
\beqn{arqpeda-model-decrease}
\barDT_{f,p}(x_k,s_k)
\geq \kappa_{\Delta m} \, \bigmin_{j\in\ii{q}}\epsilon_j^{\pi_j},
\eeqn
for every iteration of algorithm \al{AR$qp$EDA2} before termination, where 
\[
\kappa_{\Delta m} \eqdef
 \left\{\begin{array}{ll}
 \bigfrac{\sigma_{\min}}{(p+1)!}
 \left(\bigfrac{\varsigma(1-\theta)(1-\omega)}
  {2q!(L_{f,p}+\sigma_{\max})(1+\omega)}\right)^{\frac{p+1}{p-q+1}}
  & \tim{if} q\in\{1,2\}, \\*[2ex]
 \bigfrac{\sigma_{\min}}{(p+1)!}
 \left(\bigfrac{\varsigma(1-\theta)(1-\omega)\kappa_{\delta,\min}^{q-1}}
  {2q!(L_{f,p}+\sigma_{\max})(1+\omega)}\right)^{\frac{q(p+1)}{p}}
  & \tim{if} q>2,
  \end{array}\right.
\]
whith $\kappa_{\delta,\min} \eqdef \kappa_\delta(\sigma_{\max})$,
$\sigma_{\max}$ by \req{extDA-sigmaupper}, and
\[
\pi_j \eqdef
 \left\{\begin{array}{ll}
   \bigfrac{p+1}{p-j+1} & \tim{if} q\in\{1,2\},\\*[2ex]
   \bigfrac{j(p+1)}{p}  & \tim{if} q>2.\\*[2ex]
  \end{array}\right.
\]
If, in addition, $k\in\calS$,
\beqn{arqpeda-obj-decrease}
f(x_k) - f(x_{k+1}) \geq 
  (\eta_1-2\omega)\kappa_{\Delta m} \, \bigmin_{j\in\ii{q}}\epsilon_j^{\pi_j}.
\eeqn
}

\proof{The bound \req{extDA-Dphi}, \req{extDA-sigupdate} and
  Lemma~\ref{arqpeda-longs-l} together imply that for every $k\in \calS$,
  there exists a $j\in \ii{q}$ such that
\[
\barDT_{f,p}(x_k,s_k)
\geq \left\{\begin{array}{l}
\!\!\!\bigfrac{\sigma_{\min}}{(p+1)!}
\left(\bigfrac{\varsigma(1-\theta)(1-\omega)}
    {2j!(L_{f,p}+\sigma_{\max})(1+\omega)}\right)^{\frac{p+1}{p-j+1}}
   \!\epsilon_j^\frac{p+1}{p-j+1}\\*[2ex]
   \hspace*{70mm}\mbox{if } q\!\in\!\{1,2\},\\*[2ex]
\!\!\!\bigfrac{\sigma_{\min}}{(p+1)!}
\left(\bigfrac{\varsigma(1-\theta)(1-\omega)\kappa_{\delta,\min}^{q-1}}
   {2j!(L_{f,p}+\sigma_{\max})(1+\omega)}\right)^{\frac{j(p+1)}{p}}
  \!\epsilon_j^{\frac{j(p+1)}{p}}\\*[2ex]
  \hspace*{70mm}\mbox{if } q\!>\!2,
\end{array}\right.
\]
and \req{arqpeda-model-decrease} follows. Suppose now that $k$ is the index of
a successful iteration before termination.  Then, using \req{extDA-Df+-DT} and
\req{extDA-Df-DT}, the implication that $\rho_k \geq \eta_1$ with \req{extDA-rhok-def}, and 
  \req{extDA-Dphi} and \req{extDA-sigupdate},
  \[
  \begin{array}{lcl}
  f(x_k)-f(x_{k+1})
  & \geq  & [\barf(x_k) - \barf(x_{k+1})]
  - 2\omega \barDT_{f,p}(x_k,s_k) \\
  & \geq & (\eta_1-2\omega)\barDT_{f,p}(x_k,s_k)\\
  & \geq & (\eta_1-2\omega) \kappa_{\Delta m} \,\bigmin_{j\in\ii{q}}\epsilon_j^{\pi_j},
  \end{array}
  \]
  where we note that $\eta_1-2\omega > 0$ from \req{extDA-acc-init}.  This
  proves \req{arqpeda-obj-decrease}.
}

We are now in position to state formally a bound on the evaluation complexity
of the \al{AR$qp$EDA2} algorithm.

\lbthm{arqpeda-complexity}{
Suppose that \textbf{f.Bb} and \textbf{f.D0pL} hold. Suppose moreover
that the \al{AR$qp$EDA2} algorithm chooses $\delta_k$ for each $k$ 
so that the conclusions of Lemma~\ref{arqpeda-newkappadelta} hold.
\begin{enumerate}
  \item If $q\in \{1,2\}$, then
     there exist positive constants 
     $\kats{ARqpEDA2,1}$, $\kata{ARqpEDA2,1}$, $\katc{ARqpEDA2,1}$ $\kate{ARqpEDA2,1}$
     and $\katf{ARqpEDA2,1}$ such that, for any $\epsilon \in (0,1]^q$, the
     \al{AR$qp$EDA2} algorithm requires at most 
     \beqn{arqpeda-fcomp1}
     \begin{array}{lcl}
     \evbndf{ARqpEDA2,1} & \eqdef
      & \kata{ARqpEDA2,1}\bigfrac{f(x_0)-\flow}{\bigmin_{j\in\ii{q}}\epsilon_j^{\frac{p+1}{p-j+1}}}
      +\katc{ARqpEDA2,1}\\*[2ex]
     & = & \calO\left(\bigmax_{j\in\ii{q}}\epsilon_j^{-\frac{p+1}{p-j+1}}\right)
     \end{array}
     \eeqn
     evaluations of $f$, and at most
     \beqn{arqpeda-dcomp1}
     \begin{array}{lcl}
     \evbndd{ARqpEDA2,1} & \eqdef
     & \kats{ARqpEDA2,1}\bigfrac{f(x_0)-\flow}
           {\bigmin_{j\in\ii{q}}\epsilon_j^{\frac{p+1}{p-j+1}}}
      +\kate{ARqpEDA2,1} \left|\log\left(\bigmin_{j\in\ii{q}}\epsilon_j\right)\right|
      +\katf{ARqpEDA2,1}\\*[2ex]
     & = & \calO\left(\bigmax_{j\in\ii{q}}\epsilon_j^{-\frac{p+1}{p-j+1}}\right)
     \end{array}
     \eeqn
     evaluations of the derivatives of $f$ of orders one to $p$ to
     produce an iterate $x_\epsilon$ for which $\phi_{f,j}^1(x_\epsilon)\leq
     \epsilon_j/j!$ for  all $j\in\ii{q}$.
\end{enumerate}
}{
\begin{enumerate}
\setcounter{enumi}{1}
  \item If $q > 2$, then there exist positive constants $\kats{ARqpEDA2,2}$,
    $\kata{ARqpEDA2,2}$,  $\katc{ARqpEDA2,2}$,  $\kate{ARqpEDA2,2}$ and
    $\katf{ARqpEDA2,2}$ such that, for any $\epsilon \in (0,1]^q$, the
    \al{AR$qp$EDA2} algorithm requires at most 
     \[
     \begin{array}{lcl}
     \evbndf{ARqEDAp,2} & \eqdef
      & \kata{ARqpEDA2,2}\bigfrac{f(x_0)-\flow}
           {\bigmin_{j\in\ii{q}}\epsilon_j^{\frac{j(p+1)}{p}}}
     +\katc{ARqpEDA2,2} \\*[2ex]
     & = & \eo{\bigmax_{j\in\ii{q}}\epsilon_j^{-\frac{j(p+1)}{p}}}
     \end{array}
     \]
     evaluations of $f$, and at most
     \[
     \begin{array}{lcl}
     \evbndd{ARqEDAp,2} & \eqdef
      & \kats{ARqpEDA2,2} \bigfrac{f(x_0)-\flow}
           {\bigmin_{j\in\ii{q}}\epsilon_j^{\frac{j(p+1)}{p}}}
      +\kate{ARqpEDA2,2} \left|\log\left(\bigmin_{j\in\ii{q}}\epsilon_j\right)\right|
      +\katf{ARqpEDA2,2}\\*[2ex]
     & = & \calO\left(\bigmax_{j\in\ii{q}}\epsilon_j^{-\frac{j(p+1)}{p}}\right)
     \end{array}
     \]
     evaluations of the derivatives of $f$ of orders one to $p$ to 
     produce an iterate $x_\epsilon$ for which
     $\phi_{f,j}^{\delta_{\epsilon,j}}(x_\epsilon)\leq \epsilon_j \,\delta_{\epsilon,j}^j/j!$
     for some $\delta_\epsilon \in (0,1]^q$ and all $j\in\ii{q}$.
  \end{enumerate}
}

\proof{Consider first the case where  $q\in \{1,2\}$.
Using \req{arqpeda-obj-decrease} in Lemma~\ref{lem-arqpeda-model-decrease}, 
\textbf{f.Bb} and the fact that the sequence $\{f(x_k)\}$ is non-increasing, 
we deduce that the algorithm needs at most
  \beqn{arqpeda-nsucc1}
  \kats{ARqpEDA2,1}\frac{f(x_0)-\flow}{\epsmin^{\frac{p+1}{p-q+1}}}+1
  \eeqn
  \emph{successful iterations} to produce a point $x_\epsilon$ for which
  $\phi_{f,j}^1(x_\epsilon)\leq \epsilon_j/j!$ for  all $j\in\ii{q}$, where
  \beqn{arqpeda-kats1}
  \kats{ARqpEDA2,1} \eqdef
  \frac{(p+1)!}{(\eta_1-2\omega)\sigma_{\min}}
  \left( \frac{2q!(L_{f,p}+\aacc_{\max}+\sigma_{\max})(1+\omega)}{(1-\theta)(1-\omega)}\right).
  \eeqn
  We may then invoke Lemma~\ref{arqpeda-SvsU} to deduce that the total
  number of iterations required is bounded by
  \[
  |\calS_k| \left(1+\frac{|\log\gamma_1|}{\log\gamma_2}\right)
  + \frac{1}{\log\gamma_2}
  \log\left(\frac{\sigma_{\max}}{\sigma_0}\right)+1,
  \]
  where $\sigma_{\max}$ is given by \req{extDA-sigmaupper},
  and hence the total number of approximate function evaluations is at most
  twice this number, which yields \req{arqpeda-fcomp1} with the coefficients
  \[
  \kata{ARqpEDA2,1} =
  2\kats{ARqpEDA2,1}\left(1+\frac{|\log\gamma_1|}{\log\gamma_2}\right)
  \tim{ and }
  \katc{ARqpEDA2,1} = \frac{2}{\log\gamma_2}
  \log\left(\frac{\sigma_{\max}}{\sigma_0}\right)+2.
  \]
  In order to derive an upper bound on the the number of derivative
  evaluations, we now have to count the number of additional
  evaluations caused by the need to approximate the derivatives to the desired
  accuracy. Again, repeated evaluations at a given iterate $x_k$ are only
  needed when the current values  of the required absolute errors are smaller 
  than used previously at $x_k$. Recall that these required absolute errors are 
  initialised in Step~0 of the \al{AR$qp$DA} algorithm, and by construction
  decrease linearly with rate $\gamma_\aacc$ on every pass to Step~5.
  We may now use
  Lemmas~\ref{arqpeda-step1},
  \ref{arqpeda-step2}, \ref{arqpeda-newallkappadelta} and
  \ref{extDA-sigmaupper-lemma} 
  to deduce that the maximal accuracy bound $\max_{i\in\ii{p}}\aacc_{i,\countacc}$ 
  will not be reduced below
  \[
  \begin{array}{l}
  \kap{acc} \,\epsmin^{q+1}
  \eqdef \min\left[
  \bigfrac{\varsigma\omega}{4q!} 
  [\kappa_\delta(\sigma_{\max})]^{q-1},
  \kapsteptwo(\sigma_{\max}) \right]\,\epsmin^{q+1} \\*[2ex]
  \leq \min\left[
  \,\bigfrac{\omega}{4} 
  \min\left[ \bigfrac{\varsigma\epsilon_j}{8(1+\omega)
 \max[\barL_f, \sigma_k]},\delta_{k,j}^{(0)} 
     \right]^{j-1} \bigfrac{\epsilon_j}{j!},\,
     \kapsteptwo(\sigma_{k,\max}) \epsmin^{q+1} \,\right] 
  \end{array}
  \]
  at iteration $k$. As $\kap{acc}$ is independent of $k$, it follows that no
  further evaluations of $\{\overline{\nabla_x^if}(x_k)\}_{i=1}^p$ can
  possibly be required during iteration $k$ or beyond once the largest initial
  absolute error $\max_{j\in\ii{p}}\aacc_{j,0}$ has been reduced by successive
  multiplications by $\gamma_\aacc$ in Step~5 sufficiently often to ensure that
  \beqn{arqpeda-adamax}
  \gamma_\aacc^{\countacc} [\max_{i\in\ii{p}}\aacc_{i,0}]
  \leq \kap{acc}\,\epsmin^{q+1}.
  \eeqn
  Since the $\aacc_{i,0}$ are initialised in the 
  algorithm so that
  $\max_{i\in\ii{p}}\aacc_{i,0}\leq \aacc_{\max}$, the bound \req{arqpeda-adamax}
  is achieved once $\countacc$, the number of decreases in 
  $\{\aacc_{i,\countacc}\}_{j=1}^p$, is large enough to guarantee that
  \beqn{arqpeda-ieps-bound}
  \gamma_\aacc^{\countacc} \aacc_{\max} \leq  \kap{acc} \epsmin^{q+1}.
  \eeqn
  Thus we obtain that the number of evaluations of the derivatives 
  of the objective function that occur during the course of the \al{AR$qp$EDA2}
  algorithm is at most
  \[
  |\calS_k| + |\calA_k| = |\calS_k| + k_{\aacc,\min},
  \]
  i.e., the number successful iterations in \req{arqpeda-nsucc1} plus
  \[
  \begin{array}{lcl}
  k_{\aacc,\min}  & \!\! \eqdef
 \!\! & \left\lfloor
  \bigfrac{1}{\log(\gamma_\aacc)}
  \left\{(q+1)\log\left(\epsmin \right)
    +\log\left(\frac{\kap{acc}}{\aacc_{\max}}\right)\right\}
  \right\rfloor \\*[2ex]
  & \!\!\leq \!\! &
  \bigfrac{q+1}{|\log(\gamma_\aacc)|}\left|\log\left(\epsmin\right)\right|
    + \bigfrac{1}{|\log(\gamma_\aacc)|}
      \left|\log\left(\frac{\kap{acc}}{\aacc_{\max}}\right)\right|+1,
  \end{array}
  \]
  the smallest value of $\countacc$ that ensures \req{arqpeda-ieps-bound}.
  This leads to the desired evaluation bound \req{arqpeda-dcomp1}
  with the coefficients
  \[
  \kate{ARqpEDA2,1} \eqdef
  \frac{q+1}{|\log\gamma_\aacc|}
  \tim{ and }
  \katf{ARqpEDA2,1} \eqdef
  \bigfrac{1}{|\log(\gamma_\aacc)|}
    \left|\log\left(\frac{\kap{acc}}{\aacc_{\max}}\right)\right|+2.
  \]

  The reasoning is essentially the same for the case where $q>2$, except that, 
  in view  of \req{extDA-longs}, we use
  \[
  \kats{ARqpEDA2,2}
  =\bigfrac{(p+1)!}{(\eta_1 -2\omega) \sigma_{\min}}
  \left(\bigfrac{2j!(L_{f,p}+\sigma_{\max})(1+\omega)}
       {(1-\theta)(1-\omega)\kappa_{\delta,\min}^{j-1}}\right)^{\frac{p+1}{p}}
  \]
  instead of \req{arqpeda-kats1}.  This then yields 
  \[
  \kata{ARqpEDA2,2}  \eqdef
  \kats{ARqpEDA2,2}
  \left(1 +\frac{|\log \gamma_1|}{\log\gamma_2}\right),
  \]
  all other constants being unchanged, that is
  \[
  \katc{ARqpEDA2,2} = \katc{ARqpEDA2,1},
  \ms
  \kate{ARqpEDA2,2} = \kate{ARqpEDA2,1}
  \tim{ and }
  \katf{ARqpEDA2,2} = \katf{ARqpEDA2,1}.
  \]
} 

\noindent
Since the orders in $\epsilon_{\min}$ are the same as those derived for 
the \al{AR$qp$} algorithm using exact evaluations (as defined in
\cite[Chapter~12]{CartGoulToin21}) and because these were proved to be sharp 
(see Section~12.2.2.4 in this reference), the same conclusion obviously holds for
the \al{AR$qp$EDA2} algorithm using inexact evaluations.

\numsection{Conclusions}

Given the significant complexity of the theory advanced above, the reader
will undoubtly understand why a simpler version of the \al{AR$qp$EDA} has been
developed and analyzed in \cite[Chapter 13]{CartGoulToin21}.  However,
\al{AR$qp$EDA2} is not without merits.  In particular, its distinguishing
feature, the requirement \req{eda-barDT-acc}, may be of interest as it is 
independent of variable scaling, a sometimes very desirable property.
  
\end{document}